\documentclass[11pt]{article}  

\usepackage{amsmath,amsthm,amssymb}

\usepackage{color}
\usepackage{graphicx}
\usepackage{subfigure}
\usepackage{psfrag}

\usepackage{verbatim}
\usepackage{floatrow}

\def\titlerunning#1{\gdef\titrun{#1}}
\makeatletter
\def\author#1{\gdef\autrun{\def\and{\unskip, }#1}\gdef\@author{#1}}
\def\address#1{{\def\and{\\\hspace*{18pt}}\renewcommand{\thefootnote}{}%
\footnote {#1}}}%
\makeatother
\def\email#1{e-mail: #1}

\def\keywords#1{\par\medskip
\noindent\textbf{Keywords.} #1}

\def\a{{\bf a}}

\def\B{{\bf B}}

\def\j{{\bf j}}

\def\n{{\bf n}}

\def\P{{\bf P}}

\def\r{{\bf r}}

\def\t{{\bf t}}

\def\v{{\bf v}}

\def\z{{\bf z}}

\def\rd{{\rm d}}

\def\RR{{\mathop{{\rm I}\kern-.2em{\rm R}}\nolimits}}

\def\EE{{\mathop{{\rm I}\kern-.2em{\rm E}}\nolimits}}

\newcommand{\be}{\begin{equation}}
\newcommand{\ee}{\end{equation}}
\newcommand{\ba}{\begin{eqnarray}}
\newcommand{\ea}{\end{eqnarray}}
\newcommand{\bi}{\begin{itemize}}
\newcommand{\ei}{\end{itemize}}

\newtheorem{rmk}{Remark}

\begin{document} 
\titlerunning{}
\title{$C^2$ continuous  time--dependent feedrate scheduling \\ with configurable kinematic constraints} 

\author{Carlotta Giannelli, Duccio Mugnaini, Alessandra Sestini}

\date{}

\maketitle

\address{
Carlotta Giannelli, Alessandra Sestini,
Dipartimento di Matematica e Informatica ``U. Dini,'' 
Universit\`a 
di Firenze, Viale Morgagni 67/A, I--50134 Firenze, Italy,
\email{carlotta.giannelli@unifi.it, alessandra.sestini@unifi.it}
\and
Duccio Mugnaini,
Dipartimento di Scienze e Alta Tecnologia,
Universit\`{a} degli Studi dell'Insubria,
Via Valleggio 11, I-22100 Como, Italy,
\email{dmugnaini@uninsubria.it}
}



\maketitle

\begin{abstract}
We present a configurable trajectory planning strategy on planar paths for offline definition of  time-dependent $C^2$ piecewise quintic feedrates.  
The more conservative formulation ensures chord tolerance, as well as prescribed bounds on velocity, acceleration and jerk Cartesian components. Since the less restrictive formulations of our strategy can usually still ensure all the desired bounds while simultaneously producing faster motions, the configurability feature is useful not only when reduced motion control is desired but also when full kinematic control has to be guaranteed. Our approach can be applied to any planar path with a piecewise sufficiently smooth parametric representation.  
When Pythagorean-hodograph spline curves are considered, the corresponding accurate and efficient CNC interpolator algorithms can be exploited.  \end{abstract}


\section{Introduction}
Trajectory planning is a fundamental issue in several fields related to automation, such as robotics, digital animation and computer aided manufacturing. It requires the identification of a feasible path and the definition of an admissible time law. Even if conventional controllers were restricted to piecewise linear or circular paths, modern software systems enable the design of smooth curvilinear paths where zero velocity values can be avoided. Once the path is identified, the design of a suitable time law is crucial to achieve high-speed control. This task corresponds to specify the object speed along the path, generally called \emph{feedrate} function.

The final output of a trajectory planning scheme is a set of \emph{reference points} for the  control system. They indicate the prescribed machine positions with respect to uniform time spacing and constitute the input values for the controller to guide the movement of the machine. The feedrate scheduler provides  the velocity  profile for the planned path by taking into account the given geometric path and kinematic constraints. A characterizing feature of any feedrate scheduler relies on the level of kinematic control it ensures (or tries to ensure). This may strongly vary between different approaches. For example, full jerk control is usually not guaranteed, since the contribution of the jerk centripetal component is often ignored. The continuity properties of the motion represent another important aspect of the control strategy. In particular, since the time law should avoid unnecessary vibrations (that may also cause mechanical damages), smooth motions are strongly preferable. Only schedulers providing $C^2$ feedrate functions allow the possibility of ensuring acceleration and jerk continuity on sufficiently smooth path segments. 
This is an important feature of the scheduler since sudden changes or reversal in the kinematic profiles are avoided in the motions it produces, see e.g. \cite{gasparetto12}.

In view of their compatibility with standard spline representation of Computer Aided Design systems, different kind of Non-Uniform Rational B-splines (NURBS) interpolators were recently proposed. The scheduling computation can be executed before or during the interpolation phase, leading to \emph{offline} or \emph{online} approaches, respectively. While the offline interpolator avoids real-time interaction at the price of using additional storage of computed data, the online solution requires more computational power at running time. Examples of  offline and online NURBS interpolators with chord error, acceleration and jerk limitations were presented in \cite{lee11,xinhua16} and \cite{annoni12}, respectively. 

The feedrate profile is often expressed as a time-dependent function, but representations with respect to the curve parameter or the curvature are also considered.
For example, when its representation with respect to the curvature profile is considered, the chord error control is highly facilitated. This is an important feature in CNC machining and processing, see, e.g., \cite{farouki96b,yeh02}. Real-time interpolator algorithms for different feedrate variations (parameter-dependent, curvature-dependent and hybrid formulations) were also recently considered in connection with a corner rounding approach for high-speed machining \cite{farouki16b}. A time-dependent feedrate instead simplifies the development of controllers that take into account different kinematic constraints. Among the wide set of references in this setting, trigonometric feedrate profiles have been considered by several authors in connection with different continuity requirements, see e.g., \cite{lee11} for the $C^1$ case and \cite{xinhua16} for the case of higher continuity.

In this paper we present a configurable offline trajectory planning strategy that provides a time-dependent $C^2$ quintic spline feedrate profile. The more conservative formulation of the scheduler ensures chord tolerance control, as well as certified velocity, acceleration and jerk bounds. Since the less restrictive formulations of our strategy can usually still ensure all the desired bounds while simultaneously producing faster motions, the configurability feature is useful not only when reduced motion control is desired, but also when full kinematic control has to be guaranteed. In all formulations the scheduler first scans the path to detect a set of special points that subdivide the curve into successive curve blocks. The feedrate function is defined with a typical trapezoidal profile in any curve block. The maximal feedrate value of each curve block and the feedrate values at certain special points are properly initialized by taking into account the kinematic constraints and the block length. The feedrate is then defined as a time-dependent $C^2$  piecewise quintic polynomial function.
Specific bounds on the first and second derivatives of the feedrate are suitably considered in the minimization of the block traversal time. The configurable feature of the scheduler allows us to specify if only acceleration control is required, or if either tangential or full jerk limitations should be also taken into account. Another possibility offered by our strategy enables the choice between maximal and mean values of curvature in each curve block, a key geometric aspect for the design of versatile kinematic control.

In order to exploit accurate and reliable real-time interpolators, the class of Py\-tha\-gorean-hodograh (PH) curves is here considered \cite{farouki2008} for specifying the considered planar paths. The polynomial nature of the parametric speed of a PH curve enables an exact computation of its arc-length and leads  to robust processing algorithms. By avoiding repeated numerical approximations, the precision of the method is strongly enforced. The results obtained with Pythagorean-hodograh spline curves in CNC applications have confirmed the potential of this approach, see e.g., \cite{sir05a,sir07b,tsai01}, as well as  \cite{farouki2008} and references therein. The complete trajectory planning scheme should properly integrate the feedrate scheduling with the geometric part of the motion. Planar algorithms for approximating and interpolating arbitrary input data with PH spline schemes with different order of continuity are available, as for example \cite{farouki01a,fms2003}. More recently, tangent and curvature continuous path planning with obstacle avoidance techniques based on planar PH quintic spline curves were also proposed in \cite{gms15} and \cite{dgms17}, respectively. These path planning schemes can be suitably combined with the feedrate scheduler here proposed to properly define an optimal trajectory for a given path \cite{mugnaini2017}.
 
The structure of the paper is as follows. Section~\ref{sec:pre} introduces some preliminary concepts and the kinematic constraints here considered. The configurable feedrate scheduler is presented in Section~\ref{sec:sch}, while the CNC interpolator algorithm for PH curves is reviewed in Section~\ref{sec:int}. Finally, Section~\ref{sec:res} illustrates the performance of our approach on a selection of significative examples, while Section~\ref{sec:clo} concludes the paper with some final remarks.

\section{Preliminaries and kinematic constraints}\label{sec:pre}
We assume that a planar path is given in parametric form, $\r(\xi), \xi \in [\xi_L\,,\,\xi_R]$ with $\r$ sufficiently smooth. We  are interested in determining a suitable time low for  traversing the given path. More precisely, if $t$ denotes the time and  $s$ denotes the {\it arc length}, we define the {\it feedrate} $v(t)$ as 
$$v(t) \,:=\, \frac{ds}{dt} (t).$$
The corresponding global traversal time $T$ is then implicitly given by
\be \label{time_impl} 
\int_0^T v(t) \,dt \,\,=\, S\,,
\ee
where $S$ is the total length of the path. Note that, when the feedrate is given as a function of the parameter $\xi$, the global traversal time $T$ can be explicitly obtained from the following formula:
$$T = \int_{\xi_L}^{\xi_R} \frac{\sigma}{v} \, d\xi\,,$$
where $\sigma := {ds}/{d\xi} = |\r'(\xi)|$ denotes the curve {\it parametric speed}.  Since we want to prescribe a time-dependent feedrate profile, we will use the relation in (\ref{time_impl}).

Let $\v := \dot{\r}$ be the velocity vector along the trajectory (with the dot symbol denoting a derivative with respect to time). Being defined at each regular point as
\be \label{velocity} \v  = v \, \t\,, \ee
where $\t := \r'/|\r'|$ is the unit tangent vector to the curve (with the prime symbol denoting derivative with respect to the parameter $\xi$), the absolute value  $|v|$  of the feedrate is simply the velocity module.
Consequently, the absolute value of both the Cartesian components of the velocity can be bounded from above by $V_m$ by requiring 
\be \label{boundsv} 0 \le v  \le V_m\,, \ee
where $V_m$ is called the {\it commanded feedrate}. The required positivity for $v$ means that we are interested only in motions that avoid any inversion of direction.   

The acceleration vector $\a := \ddot{\r}$ is well defined at any regular point and its expression is given by
\be \label{acceleration} \a = \dot{v} \t \,-\, v^2 \kappa\, \n\,, \ee
since ${d\,\t}/{ds} = -\kappa\,\n,$ where  
$$
\n := \t \times \z \quad \mbox{(unit normal)}\quad 
\mbox{and } \quad \kappa := \frac{(\r' \times \r'') \cdot \z}{\sigma^3} \quad \mbox{(curvature)}$$
with $\z$ denoting a unit vector orthogonal to the plane of the path. 

By differentiating the expression of $\a$ given in (\ref{acceleration}), we obtain the expression of the {\it jerk} vector $\j := \dddot{\r}$  at any regular point
\be \label{jerk}  \j \,=\, (\ddot{v} -v^3\kappa^2) \t \,-\, (3v\dot{v}\kappa + v^3 \kappa' /\sigma) \n\,,\ee
in terms of the feedrate, its first and second derivative, and other geometric quantities connected to the given path.

Since in computer numerical control the motion is separately prescribed on each axis, the general interest is in bounding the absolute value of  the Cartesian components of the acceleration vector $\a$ and, possibly, also those of the jerk vector $\j.$ In view of (\ref{acceleration}), we can guarantee that both $|a_x| \le A_m$ and $|a_y| \le A_m\,$ only by requiring
\begin{equation} \label{acccont} 
|\dot v|\le A_t : = p_a \, A_m\,, \qquad v^2 |\kappa| \le A_c := \sqrt{1-p_a^2}\,A_m,
\end{equation} 
with $p_a \in (0\,,\,1)$. A typical choice is $p_a = 1/\sqrt{2}.$  Analogously, from the jerk expression in (\ref{jerk}) we obtain that  $|j_x| \le J_m$ and $|J_y| \le J_m\,$ is guaranteed only if
\begin{equation} \label{jerktcont} 
|\ddot{v}| \le J_{t_1} \,: = q_j p_j \, J_m\,, \qquad
 v^3 \kappa^2 \,\le J_{t_2} \,:= (1-q_j)\,p_j J_m\,, 
\end{equation} 
and
\begin{equation} \label{jerkccont} 
  3|\kappa| A_t v + |\kappa'| / \sigma v^3  \le J_c := \sqrt{1-p_j^2} J_m\,,
\end{equation} 
with  $p_j, q_j \in (0\,,\,1)$. Typical choices in this case are $p_j = 1/\sqrt{2}$ and $q_j = 1/2.$ Note that the inequality in (\ref{jerkccont}) related to the jerk centripetal component control, is often not considered in the literature, see  \cite{annoni12} for one example. As noticed in the appendix of \cite{annoni12},  its left hand side  is an increasing polynomial in $v.$ Thus there exists exactly
one positive value $v_{jc}$ for $v$ such that  (\ref{jerkccont}) holds if and only if $v \le v_{jc}\,.$  When $\kappa' \ne 0,$ the expression of $v_{jc}$ is given by Cardano's rule,
\begin{equation} \label{card} v_{jc} = v_{jc}(\kappa,w) := \sqrt[3]{\lambda_+} +  \sqrt[3]{\lambda_-}\,, \quad \mbox{with } \lambda_{\pm} = \frac{J_c }{|w|} \pm \sqrt{\left(\frac{J_c }{2 |w|} \right)^2 + 
\left(\frac{\kappa A_t }{ |w|} \right)^3}\,,
\end{equation} 
where $w := \kappa'/\sigma$ is the curvature derivative with respect to the arc length.

\begin{figure}[h!]
\begin{center}
\includegraphics[scale=0.6, trim = 0cm 1.5cm 0cm 1.5cm, clip]{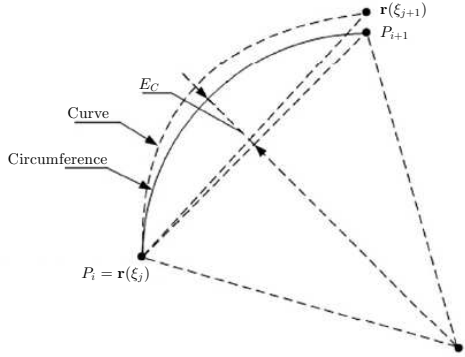}
\caption{Chordal error $E_c(j)$ at the interpolated point $\r(\xi_j)$.}
\label{fig:chordal_error}
\end{center}
\end{figure}

In addition to the kinematic control,  it is also natural to ensure a suitable control on the {\it chord error} $E_c$ associated with a given {\it sampling time} $\Delta t$  during the overall traversal of the path. Since the chord error is a discrete function whose $j$-th value $E_c(j)$  is associated with  $t_j =j \Delta t$, it necessarily depends on the sampling time $\Delta t.$ As shown in Figure \ref{fig:chordal_error}, when a point with non vanishing curvature is considered, $E_c(j)$ is defined as  
 the distance between the osculating circumference through the corresponding {\it interpolated point} ${\bf P}_j = \r(\xi_j)$   and a related circumference chord,
  \be \label{chordal} E_c(j) \,:=\,\left\{ \begin{array}{lll } \frac{1}{|\kappa_j|} - \sqrt{\frac{1}{\kappa_j^2} - \frac{L_j^2}{4}}\,, & \mbox{if } \kappa_j \ne 0\,, \cr
  0\,, & \mbox{otherwise,} \end{array} \right.
 \ee
with   $L_j  \,:=\,  v(t_j) \Delta t$. This means that $L_j =  \| {\bf P}_{j+1} - {\bf P}_j \|_2\,, $ where ${\bf P}_{j+1}$ is the next interpolated point along the circumference, by assuming a constant feedrate equals to $v(t_j)$ and considering that the tool moves on a straight line between any pair of successively interpolated points. Observe that, if $\kappa_j \ne 0,$ $E_c(j)$ is well defined only if $L_j \le 2 / |\kappa_j|,$ that is only if $v(t_j) |\kappa_j| \Delta t \le 2.$ 

A suitable requirement is
\be \label{chcontr} E_c(j) \le D\,,\ee
where $D$ is an assigned positive threshold. Clearly, if $D \ge 1/|\kappa_j|$ then (\ref{chcontr}) holds. On the other hand, when   $D< 1/|\kappa_j|,$ it can be verified with some computations that the inequality in (\ref{chcontr}) corresponds to require the following condition on the feedrate
\be \label{condchord} v(t_j) \le   \frac{2}{\Delta t}\, \sqrt{\frac{1}{\kappa_j^2}-\left(\frac{1}{|\kappa_j|} - D \right)^2}\,. \ee 
It is interesting to note also that, if $v(t_j) < {2}/{\Delta t |\kappa_j|}\,,$ it is reasonable to use  a first order Taylor expansion for the square root in (\ref{chordal}) with respect to $\Delta t^2,$ leading to
$$ E_c(j) \,\approx\, \frac{1}{|\kappa_j|} \left[1-\frac{1}{2} \left( 2 -   \frac{\kappa_j^2 v_j^2 \Delta t^2 }{4} \right) \right] \,=\, \frac{ \Delta t^2}{8}\ \,  |\kappa_j| v_j^2\,.$$
In this case the chord error is then approximately proportional to the absolute value of the normal component of the acceleration. 

As previously mentioned in the introduction, we present a configurable strategy that allows six alternative formulations with different level of kinematic constraints. These alternatives are labelled as $R_0, R_1, R_2$ and $S_0$, $S_1$, $S_2$ and can be collected into two groups: the \emph{relaxed} ($R_0, R_1, R_2$) and the \emph{strict} ($S_0, S_1, S_2$) formulations.
In any case, the path is preliminarily divided into segments jointly together at tangent or curvature discontinuous points, where the motion is halted. Each segment is then further subdivided into blocks by detecting on each segment a set of special points.
In all the relaxed formulations, the mean value of the absolute curvature $\kappa$ is determined on each block. In the $R_2$ formulation the mean value of $w$ is also considered. In order to obtain small block traversal times, these values are used to imposed the selected kinematic constraints in the block only with respect to a mean value. On the other hand, in all strict formulations the maximal values of the same geometric quantities are determined in each block, in order to ensure \emph{a priori} the desired kinematic constraints. These always include the bound in (\ref{condchord}) for chord error control.
In addition, the $R_0$ and $S_0$ formulations take into account only the bounds in (\ref{acccont}) which concern acceleration, while $R_2$ and $S_2$ adds to (\ref{condchord}) and (\ref{acccont}) also the bounds in (\ref{jerktcont}) and  (\ref{jerkccont}) for full jerk control. $R_1$ and $S_1$ are hybrid formulations since they add to  (\ref{condchord}) and (\ref{acccont}) only the bounds related to the tangential component of the jerk given by (\ref{jerktcont}).
  
\section{$C^2$ time--dependent feedrate scheduling}\label{sec:sch}

In any possible configuration and for all curve segments, the scheduler defines a feedrate function $v$ such that
\begin{equation} \label{vderineq}
0 \le v \le V_m\,, \quad  |\dot{v}| \le A_t\,, \quad |\ddot{v}| \le J_{t_1} \,.
\end{equation}
These inequalities, however, are not sufficient to ensure the desired chord, acceleration, and jerk control. In addition to them, by considering the expressions in (\ref{acceleration}), (\ref{jerk}), and (\ref{chordal}), we can easily derive that, to ensure full kinematic control, the feedrate $v$ should satisfy the upper bound $v \le V_b(\kappa,w),$  where $w = \kappa'/\sigma$ and
\be \label{vbound} 
V_b(\kappa,w) :=\,  
\min \left\{  \frac{2}{\Delta t}\, \sqrt{\frac{1}{\kappa^2}-\left(\frac{1}{|\kappa|} - D \right)^2}\,,\, \sqrt{\frac{A_c}{|\kappa|}}\,,\, \sqrt[3]{\frac{J_{t_2}}{\kappa^2}}\,, \,\, v_{jc}(\kappa, w)\right \}  \,,
\ee
with $v_{jc}(\kappa,w)$  defined as in (\ref{card}). Note that $\kappa'$ is assumed continuous in the segment of reference. Clearly, if we are not interested in guaranteeing jerk bounds, the minimum on the right hand side of (\ref{vbound}) is just between the first two terms. 

Let us now present the strategy of our scheduler. 
The curve is preliminarily decomposed into {\it segments} whose extrems $ \B_j, j=1,\ldots,N-1$, usually called {\it breakpoints}, are assigned points where the curve is just $G^0,$ i.e. corner points, augmented by the initial and final positions $\B_0 = \r(\xi_L)$ and $\B_N = \r(\xi_R).$  At each breakpoint we always set to zero $v, \dot v$, and also $\ddot v.$  For brevity, in what follows, we will focus on the generic $j$-th curve segment by removing  the subscript $(j)$. We will also assume that the time is reported to zero whenever a new segment is encountered.

On each curve segment the scheduler first detects a set of special points. These points allow us to subdivide the segment into blocks where a feedrate trapezoidal profile is defined. The special points are of different type for the relaxed and strict formulations and they are called {\it critical} and {\it crossing} points, respectively. In the first case, they correspond to points where the absolute curvature has a local maximum exceeding a critical curvature value $\kappa_{cr}$. In the second case instead, they correspond to points where $|\kappa| -\kappa_{cr}$ has a sign change. The critical curvature is defined as 
\be \label{critical} \kappa_{cr} :=\,  \min \left\{\frac{8 D}{V_m^2 \Delta t^2 + 4 D^2}\,, \, \frac{A_c}{V_m^2}\,,\, \sqrt{\frac{J_{t_2}}{V_m^3}} \right\}\,, 
\ee
for $R_1, S_1, R_2$ and $S_2$. By assuming (\ref{vderineq}), this is the curvature value such that, if a point of the path is encountered with a feedrate greater than $V_m,$ one among the chord error, the normal component of acceleration, and the tangential component of the jerk can exceed the related bounds. On the other hand, for the formulations $R_0$ and $S_0$, $\kappa_{cr}$ reduces to the minimum between the first two terms on the right hand side of (\ref{critical}), since they do not care of jerk bounds.  The special points are denoted in the following as $ \P_1 , \ldots, \P_{M -1}$,  while $\P_0$ and  $\P_M $ will denote the endpoints of the considered segment.  
  
\subsection{The \emph{relaxed} formulations of the scheduler}

The first configurations of the scheduler are called \emph{relaxed} (R) formulations because they consider only mean values for $\kappa$ and $w$ to try to ensure the desired kinematic control on each curve block. As mentioned in the previous section, $R_0$ takes into account only acceleration bounds, while $R_1$ adds also the two bounds in (\ref{jerktcont}) and $R_2$ requires both (\ref{jerktcont})  and  (\ref{jerkccont}).  

\subsubsection{Feedrate initialization at critical points}

We start by initializing the feedrate values at any critical point $\P_i, i=0,\ldots,M.$ At  $\P_0$ and $\P_M$   we prescribe zero values for $v$, $\dot{v}$, and $\ddot{v}.$   At any  $\P_i, i=1,\ldots,M-1,$ instead, a  sufficiently small positive value $V_i < V_m$ of the feedrate is assigned for $v,$ while $\dot{v}$  and $\ddot{v}$ are still set to zero.   In particular, in order to ensure the desired kinematic control at these points, the value $\hat V_i < V_m$ is chosen with $R_1$ and $R_2$ to initialize $V_i$
\be \label{Vifirst} \hat V_i :=\, V_{r_1}(\kappa_i)\,,
\ee
where $\kappa_i$ is the value of $\kappa$ at $\P_i$ and
\be \label{Vridotta1}
V_{r_1}(\kappa)  := \min \left\{  \frac{2}{\Delta t}\, \sqrt{\frac{1}{\kappa^2}-\left(\frac{1}{|\kappa|} - D \right)^2}\,,\, \sqrt{\frac{A_c}{|\kappa|}}\,,\, \sqrt[3]{\frac{J_{t_2}}{\kappa^2}}\,, \right \} \,.\ee
Not being interested in jerk control,  the $R_0$ formulation fixes $\hat V_i = V_{r0}(\kappa_i),$ where

\be \label{Vridotta2}
V_{r_0}(\kappa)  := \min \left\{  \frac{2}{\Delta t}\, \sqrt{\frac{1}{\kappa_i^2}-\left(\frac{1}{|\kappa|} - D \right)^2}\,,\, \sqrt{\frac{A_c}{|\kappa|}}\,  \right \}\,. \ee
Since $\kappa'$ vanishes at critical points and $\dot{v}$ is set to zero, the centripetal jerk vanishes at these points.     Furthermore, since also $\ddot{v}$ is  set to zero, the tangential acceleration vanishes. This implies that  in (\ref{Vridotta1}) we could replace $A_c$ and $J_{t_2}$ with $A_m$ and $J_m,$  respectively, see for example \cite{yeh02}. We did not consider this substitution in our scheduler since no advantages were observed in the experiments by its implementation.\footnote{In our tests we considered paths that admit a $C^2$ parametric spline representation. Some critical points were often fixed at points where $w$ is discontinuous and can remain high in the left and/or right neighborhood of these points.}

Note, however, that the $V_i$ values  can be further reduced for compatibility reasons connected to the lengths $S^{(i)}, i=1,\ldots,M\,,$ see Remark \ref{rmk:feedvar}. 

\subsubsection{Initialization of maximal feedrate value on a curve block}
\label{subsec:init_vm}
Once the needed information at the critical points are prescribed, we need to  specify the feedrate function during the motion between $\P_{i-1}$ and $\P_i$ accordingly. By following \cite{lee11}, we  require a {\it trapezoidal} feedrate profile, i.e.   $\dot{v} >0$ (acceleration phase), then $\dot{v} = 0$ (constant phase),  and finally   $\dot{v} <0$ (deceleration phase). However, we use a different analytical definition of the feedrate which produces a $C^2$ profile.   In order to reduce the traversal time, we plan the feedrate so that its central constant phase at maximal feedrate  $v = V_c^{(i)}$ is as long as possible. The maximal block feedrate value $V_c^{(i)}$  is initialized with $\hat V_c^{(i)}$ defined as follows,
\be
 \label{Vc} 
 \hat V_c^{(i)} \,:=\  \left\{ \begin{array}{ll}V_b(\kappa_c^{(i)}\,,\,w_c^{(i)}) & \mbox{in case of } R_2\,, \cr
V_{r_1}(\kappa_c^{(i)}) & \mbox{in case of } R_1 \,, \cr \
V_{r_0}(\kappa_c^{(i)}) & \mbox{in case of } R_0\,,
\end{array} \right. 
\ee
where $\kappa_c^{(i)}, w_c^{(i)}$ are the mean absolute values of $\kappa$ and $w$  between $\P_{i-1}$ and $\P_i.$  Observe that with the $R_2$ formulation it can happen that
$V_c^{(i)} < \hat V_{i-1}$  or $V_c^{(i)} < \hat V_i$. Thus, after initializing all the block maximal feedrates, we modify the initialization of the feedrate values at critical points as follows,
$$\hat V_i = \min \{\hat V_c^{(i)}, \hat V_c^{(i+1)}, \hat V_i \} \,.$$ 
Note also that, as in the case of the feedrate values at the critical points, also  $V_c^{(i)}$ can be further reduced by the algorithm for compatibility reasons connected to the length $S^{(i)}$ (see the details below).
In any case, the feedrate in the $i$-th block, $i=1,\ldots,M,$ will have  a central constant profile $v = V_c^{(i)}$ as long as possible, with $V_c^{(i)}$ as near as possible to $\hat V_c^{(i)}.$  In some cases, depending on the values of $S^{(i)}$, $V_{i-1}$, and $V_i,$  it can also happen that the three phases collapse to just two (or even one) phases. 

\subsubsection{Time-dependent feedrate scheduling on a curve block}
In order to  design a trapezoidal feedrate profile  capable to ensure the continuity of the velocity, the acceleration, and the jerk vectors at all points of $\r$ where $\kappa'$ is continuous,  we require  $v$ to be a $C^2$  quintic spline function, that can easily be represented in piecewise B\'ezier form. We then introduce our time-dependent feedrate definition by firstly associating with  the $i$-th block the time interval
$$
[ t_s^{(i)}\,,\,  t_s^{(i)} + T^{(i)}],
$$
where $ t_s^{(i+1)} := t_s^{(i)} + T^{(i)}$, with $t_s^{(1)} := 0$, and
$$
T^{(i)} := T_{acc}^{(i)} + T_{con}^{(i)} + T_{dec}^{ (i)},
$$
with  $T_{acc}^{(i)}$, $T_{con}^{(i)}$, and $T_{dec}^{(i)}$ denoting the time of the acceleration phase, the time  of the  constant maximal feedrate, and the time  of the deceleration phase, respectively.
 
Let us focus on the general $i$-th block. Clearly, the corresponding traversal time $T^{(i)}$ cannot be a priori computed, since it depends on the feedrate function $v^{(i)}(t)$, defined as 
$$ v^{(i)}(t) := \left\{\begin{array}{ll} v_{acc}^{(i)}(t) & \mbox{if } t_s^{(i)} \le t < t_s^{(i)} + T_{acc}^{(i)}\,,
 \cr 
 v_{con}^{(i)}(t) & \mbox{if } t_s^{(i)} + T_{acc}^{(i)} \le t <  t_s^{(i)} + T_{acc}^{(i)} + T_{con}^{(i)}\,,
 \cr
v_{dec}^{(i)}(t) & \mbox{if } t_s^{(i)} + T_{acc}^{(i)} + T_{con}^{(i)} \le t < 
  t_s^{(i)} + T^{(i)}   \,,
 \cr
\end{array}
\right. $$
where in the constant phase we fix $v_{con}^{(i)}(t) \equiv V_c^{(i)}.$ In particular, we assume the following definitions for the feedrate profile in the acceleration phase,     
\be \label{vaccdef} v_{acc}^{(i)} (\tau) :=\, V_{i-1} \sum_{j=0}^2 b_j^5(\tau)   \,+\,
V_c^{(i)} \sum_{j=3}^5 b_j^5(\tau) \,, \ee
 with 
$$ 
\tau \,:=\,  \frac{t-  t_s^{(i)}}{ T_{acc}^{(i)}} \,, $$
and 
$$ b_j^n(\tau) := \binom{n}{j}\, \tau^j\,(1-\tau)^{n-j},\qquad j =0,\ldots,n,$$
denoting the $j$-th Bernstein polynomial  of degree $n$.   
Analogously, in the deceleration phase, we set
  \be \label{vdecdef} v_{dec}^{(i)} (\tau) :=\, V_c^{(i)} \sum_{j=0}^2 b_j^5(\tau)   \,+\, V_i  \sum_{j=3}^5 b_j^5(\tau)  \,, \ee
with 
$$ \tau \,:=\, 
\frac{t-  t_s^{(i)}  - T_{acc}^{(i)}  - T_{con}^{(i)}}{T_{dec }^{(i)}}\,. $$ 
Let us consider the constraints we need to impose in order to ensure the fulfillment of the inequalities in \eqref{boundsv} and \eqref{acccont}, as well as of the following compatibility condition,
\be \label{comp} S_{acc}^{(i)} + S _{dec}^{(i)} \le S^{(i)}  \,, \ee
where $S_{acc}^{(i)}\, (S_{dec}^{(i)})$ denotes the length of the portion of the $i$-th block of $\r$ done in the acceleration (deceleration) phase. Note that, since $0 \le V_j< V_c^{(i)} \le V_m\,, j=i-1,i\,,$ considering the properties of the Bernstein basis we can deduce that the inequality in \eqref{boundsv} is surely satisfied on both sides. Since
$$\dot{v}_{acc}^{(i)}(t) \,=\, \frac{1}{T_{acc}^{(i)}}\, \frac{dv_{acc}^{(i)}}{d\tau} \,=\,  5\, \frac{V_c^{(i)}-V_{i-1}}{T_{acc}^{(i)}}\, b_2^4(\tau),
$$
we have $\dot{v}_{acc}^{(i)} \ge 0.$
By considering the two inequalities in \eqref{acccont}, we obtain
$$ 5\,\frac{V_c^{(i)}-V_{i-1}}{T_{acc}^{(i)}} \, \max_{\tau \in [0\,,\,1]} b_2^4(\tau) \le p_a A_m\,, \qquad
  5\, \frac{V_c^{(i)}-V_{i-1}} {(T_{acc}^{(i)})^2}\, \max_{\tau \in [0\,,\,1]}| (b_2^4)'(\tau)| \le  J_m,$$
which in turn, in view of the expression of $b_2^4(\tau)$,  reduce to
  $$ \frac{V_c^{(i)}-V_{i-1}}{T_{acc}^{(i)}} \,  \le \frac{8}{15}\, p_a A_m\,, \qquad 
  \frac{V_c^{(i)}-V_{i-1}} {(T_{acc}^{(i)})^2}\,  \le  \frac{\sqrt{3}}{10}\, J_m\,.$$
If $V_c^{(i)}$ is set equal to $x^{(i)},$ with  $\hat V_c^{(i)} \ge x^{(i)} \ge \max\{V_{i-1}\,,\,V_i\}\,,$ the value of $T_{acc}^{(i)}$ has to be fixed as 
\be \label{Taccdef} T_{acc}^{(i)} = T_{acc}^{(i)}(x^{(i)}) :=  \max \left\{ \frac{15}{8 p_a A_m} \,(x^{(i)}-V_{i-1}) \,,\,  \sqrt{ \frac{10} {\sqrt{3} J_m} \, (x^{(i)}-V_{i-1})}  \right\} \,,\ee
in order to minimize the time needed to arrive at the constant phase. With analogous computations we can verify that $\dot{v}_{dec}^{(i)} \le 0$ and $T_{dec}^{(i)}$ has to be fixed as 
\be \label{Tdecdef} T_{dec}^{(i)} = T_{dec}^{(i)} (x^{(i)}) :=  \max \left\{ \frac{15}{8 p_a A_m} \,(x^{(i)}-V_i) \,,\,  \sqrt{ \frac{10}{\sqrt{3} J_m} \, (x^{(i)}-V_i) }   \right\},\ee 
in order to ensure that the inequalities in (\ref{acccont}) hold also in the deceleration phase with  minimal time for the deceleration.
In equations (\ref{Taccdef}) and (\ref{Tdecdef}) the dependency of $T_{acc}^{(i)}$ and $T_{dec}^{(i)}$ on the value $ x^{(i)} \in [\max\{V_{i-1}\,,\,V_i\}\,,\,\hat V_c^{(i)}]$ to be selected for $V_c^{(i)}$ is emphasized.  In particular, we will fix $V_c^{(i)}$ as the maximum value of $x^{(i)}$ in this interval which guarantees the fulfillment of (\ref{comp}). 
Now, since (\ref{time_impl}) implies 
$$S _{acc}^{(i)} \,=\,T_{acc}^{(i)} \, \int_0^1 v_{acc}^{(i)}(\tau) \, d\tau\,, \qquad S _{dec}^{(i)} \,=\,T_{dec}^{(i)} \, \int_0^1 v_{dec}^{(i)}(\tau) \, d\tau\,,$$
we need to impose that
$$ T_{acc}^{(i)} \, \int_0^1 v_{acc}^{(i)}(\tau) \, d\tau\, \,+\, T_{dec}^{(i)} \, \int_0^1 v_{dec}^{(i)}(\tau) \, d\tau \,\le S^{(i)}\,, $$
which, considering (\ref{vaccdef}) and (\ref{vdecdef}), in view of the integration formulas for Bernstein polynomials is equivalent to
 \be \label{gfun}  g^{(i)}(x^{(i)}) \,\le\, 2S^{(i)}\,,\ee
with
\be \label{gidef} g^{(i)}(x^{(i)}) \,:=\, T_{acc}^{(i)}(x^{(i)}) [ V_{i-1}+ x^{(i)}] \,+\, T_{dec}(x^{(i)}) [ V_i+ x^{(i)}]\,.\ee
If $x^{(i)} = \hat V_c^{(i)}$ fulfills (\ref{gfun}), 
in order to keep the maximal admissible feedrate in the constant phase, we choose $V_c^{(i)} = \hat V_c^{(i)}$ with $\hat V_c^{(i)}$ defined in (\ref{Vc}). Otherwise, considering that $g^{(i)}$ is an increasing function with respect to $x^{(i)},$ we preliminarily check if the following inequality holds, 
\be \label{prev} g^{(i)}(\max\{V_{i-1},V_i \} ) \le 2S^{(i)}\,.\ee
If this is the case, we fix $V_c^{(i)}$ as the unique solution of the equation    $g^{(i)}(x^{(i)}) = 2S^{(i)}$ in the interval  $[\max\{V_{i-1},V_i \} \,,\, \hat V_c^{(i)}]\,.$
\begin{rmk} \label{rmk:feedvar} 
The fulfillment of the inequality in  (\ref{prev}) for any $i$-th curve block, $i=1,\ldots,M,$   has to be preliminarily checked.  If the inequality does not hold, we modify the initial feedrate values $\hat V_i, i=1,\ldots,M$, to lower values $V_i, i=1,\ldots,M,$ determined by solving the following constrained minimization problem,
$$\min \sum_{i=1}^{M-1}(V_i - \hat V_i)^2 \,\, \mbox{with } \hat V_i \ge V_i > 0\,,\,\,  g^{(i)}(\max\{V_{i-1},V_i \} )  - 2S^{(i)} \le 0.$$ 
This problem can be numerically solved by an iterative minimization procedure by initializing each $V_i$ with  $ {\min_{i=1,\ldots,M}} \hat V_i.$
 \end{rmk}

 \subsection{The \emph{strict} formulations of the scheduler} 

The formulations in the second set of configurations are called \emph{strict} (S) because they impose the selected bounds everywhere. Obviously, this is obtained at the price of defining slower motions.
 
On each segment we determine a set of inner {\it crossing points} $ \P_1 , \ldots, \P_{M -1}$, which are points where the quantity $|\kappa| - \kappa_{cr}$  has a sign change. We also denote as $\P_0$ and  $\P_M $ the end points of the considered segment. The curve portion between any pair of successive crossing points of the considered segment is again referred to as a {\it curve block} and $S^{(i)}$   denotes the length of the $i$-th curve block,   $ i=1,\ldots,M,$. The block length has to be preliminarily computed for any segment. The philosophy for the feedrate determination is somehow reversed with respect to the relaxed strategy. In this case we first determine a safe maximal feedrate for any curve block, for then selecting the feedrate values at the crossing points.
 
The maximal feedrate $V_c^{(i)}$ chosen on the $i$-th block, $i=1,\ldots M,$ is initialized  as
$$ 
\hat V_c^{(i)} 
:= \left\{ \begin{array}{ll}V_b(\kappa_{max}^{(i)}\,,\,w_{max}^{(i)}) & \mbox{in case of } S_2\,, \cr
V_{r1}(\kappa_{max}^{(i)}) & \mbox{in case of } S_1 \,, \cr \
V_{r0}(\kappa_{max}^{(i)}) & \mbox{in case of } S_0\,,
\end{array} \right. 
$$
where $\kappa_{max}^{(i)}$ and $w_{max}^{(i)}$ are the maximum values assumed by $|\kappa|$ and $|w|$ in the considered block. At the crossing point $\P_i, i=1,\ldots,M-1$, we initialize the $V_i$ with $\hat V_i$ fixed as follows, 
$$ \hat V_i := \min \{ \hat V_c^{(i)}, \hat V_c^{(i+1)} \} \,.$$
Since the inequality in  (\ref{prev}) still needs to be fulfilled, Remark \ref{rmk:feedvar} applies again.  

The feedrate inside each curve block is then selected with a trapezoidal profile by using the same approach used for the relaxed formulations. With this strategy at least one between the acceleration and deceleration phases disappears in many blocks.

\section{Application to interpolators for PH curves}\label{sec:int} 

\begin{figure}[!t]
\begin{center}
\includegraphics[scale=0.57]{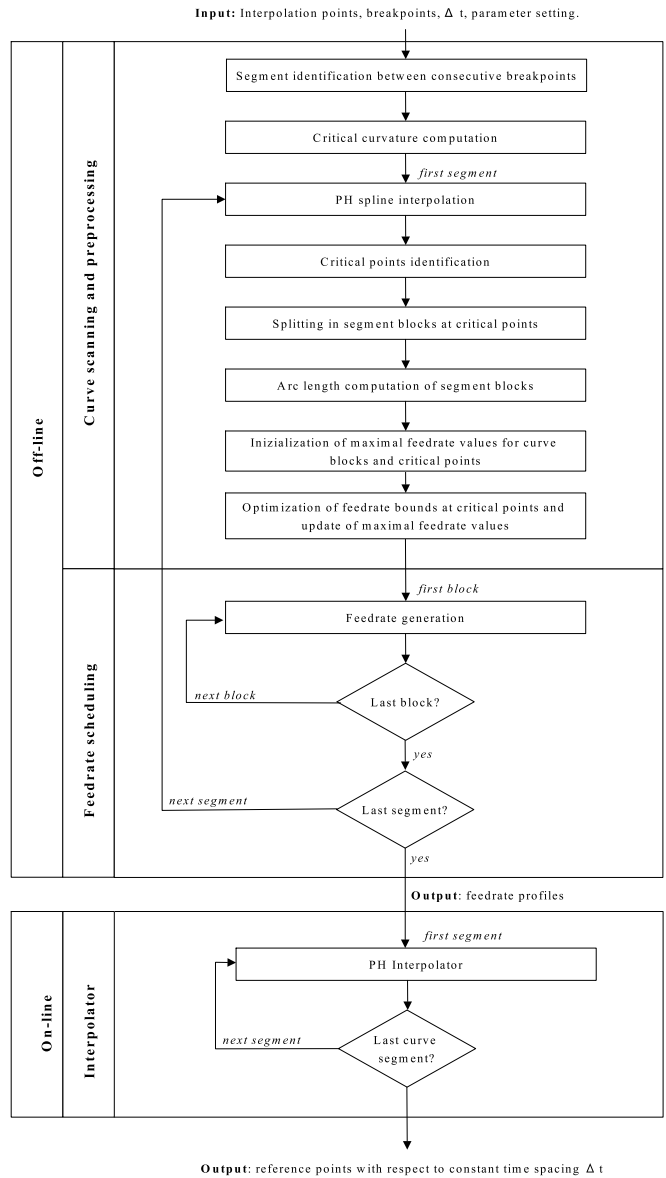}
\caption{Diagram of the relaxed strategy for feedrate definition.}
\label{fig:diagram}
\end{center}
\end{figure}

In this paper we are interested in applying the introduced feedrate scheduler on paths described by Pythagorean-hodograph curves since, among the different parametric curves commonly used in computer aided geometric design, PH curves strongly facilitates the design of efficient CNC interpolator algorithms. To avoid loss of generality,  we assume that a given planar curve is preliminarily approximated by $C^2$ quintic PH-splines. The interpolation algorithm used for this aim is the one introduced in \cite{ farouki01a}. The corresponding trajectory planning scheme is summarized in Figure \ref{fig:diagram}. For brevity, only the relaxed strategy is shown in the diagram.

In the following subsections we briefly introduce planar PH curves and the related interpolator, see \cite{farouki2008} for the details.

\subsection{Planar Pythagorean-hodograph curves}

Planar polynomial curves ${\bf r}(\xi)=(x(\xi),y(\xi)), \xi \in [0,1]$ with polynomial parametric speed $\sigma(\xi) = |\r'(\xi)|$ so that
\begin{equation}
\label{eq:planarph}
|\,{\bf r}'(\xi)\,|^2\,=\,x'^2(\xi)+y'^2(\xi)\,=\,\sigma^2(\xi)\,,
\end{equation}
are known as Pythagorean-hodograph curves \cite{farouki2008}. Although PH curves of degree $n$ have just $n+3$ degrees of freedom
 (compared to $2n+2$ for general degree-$n$ polynomial curves), 
the Pythagorean condition
(\ref{eq:planarph}) ensures distinctive
properties such as polynomial arc-length functions, rational curvature and unit tangent
and normal vectors, closed--form expression for the bending
energy, rational offset curves. For a deeper and detailed discussion see \cite{farouki2008}.

In order to satisfy condition (\ref{eq:planarph}), the three polynomials
$x'(\xi)$, $y'(\xi)$ and $\sigma(\xi)$ must comprise a Pythagorean triple. Hence,
according to the result of Kubota \cite{kubota72}, these polynomials must be
expressible in terms of other real polynomials $u(\xi)$, $v(\xi)$ and $h(\xi)$ in the
form
\begin{equation}
\label{eq:rpplanrph}
x'(\xi)\,=\,h(\xi)\left[u^2(\xi)-v^2(\xi)\right],\qquad
y'(\xi)\,=\,2\,h(\xi)\,u(\xi)\,v(\xi)\,, 
\end{equation}
with corresponding polynomial parametric speed
\begin{equation}
\label{eq:phsigma}
\sigma(\xi)\,=\,h(\xi)\left[u^2(\xi)+v^2(\xi)\right].
\end{equation}
Considering that $\sigma(\xi) = |\r'(\xi)|,$ the arc-length function is defined as
$$
s(\hat{\xi})=\int_0^{\hat{\xi}} \sigma(\xi)\,\rd\xi 
$$
which is a polynomial function of the parameter, and, consequently, it can be \textit{exactly} computed.


\subsection{PH curve interpolators}
PH interpolators based on time-dependent feedrate function
$v(t)$ are exploited to generate smooth motion trajectories for a PH spline curve $\r(\xi)$.
In order to be efficient, they require  $v(t)$ for $t \in [0, T]$ to be non-negative and to have a
simple expression $F(t)$ for its indefinite integral such that $v(t) =
\dot{F}(t)$, as in the case of our scheduler.
The PH curve interpolator implementation is based on the  relation:
\begin{equation} \label{main_relation}
s(\xi) = F(t).
\end{equation}
By considering a sampling time $\Delta t$ and a finite sequence of time
values $t_0 =\Delta t$, $t_1 = 2\Delta t$, \ldots, $t_N = N\Delta t$, we can then obtain from (\ref{main_relation}) the related sequence of parameter values, $0 = \xi_0 < \xi_1<
\ldots, <\xi_N = 1$ as the root of
$$
s(\xi_k) - F(k\Delta t) = 0 , \qquad k=1,\ldots,N.
$$
Each $\xi_k$ corresponds to a reference point $\r(\xi_k)$ on the PH curve. It will be traversed at the corresponding $t_k$.
The iterative Newton-Raphson method can be exploited to find the unique real root of the above equation, namely
\begin{equation}
\xi_k^{(r+1)} = \xi_k^{(r)} - \frac{s(\xi_k^{(r)}) - F(k\Delta t)}{\sigma
(\xi_k^{(r)})} \label{alg:interpolation}
\end{equation}
with initial guess $\xi_k^{(0)} = \xi_{k-1}$.


For the experiments presented in the next section, we implemented a PH interpolator for PH quintic
splines, in line with our choice to assign the planar paths. Particular care is needed whenever the $r$-th
approximation of $\xi_k$ in (\ref{alg:interpolation}) results greater than $1$. This corresponds to a change of polynomial  segment in the PH spline at $\xi_k$. Few iteration of the 
The Newton-Raphson method was always sufficient to obtain approximations of the reference points of the desired accuracy.



\section{Numerical results}\label{sec:res}

\begin{figure}[!th]
\begin{center}
\subfigure[]{\includegraphics[scale=0.30]{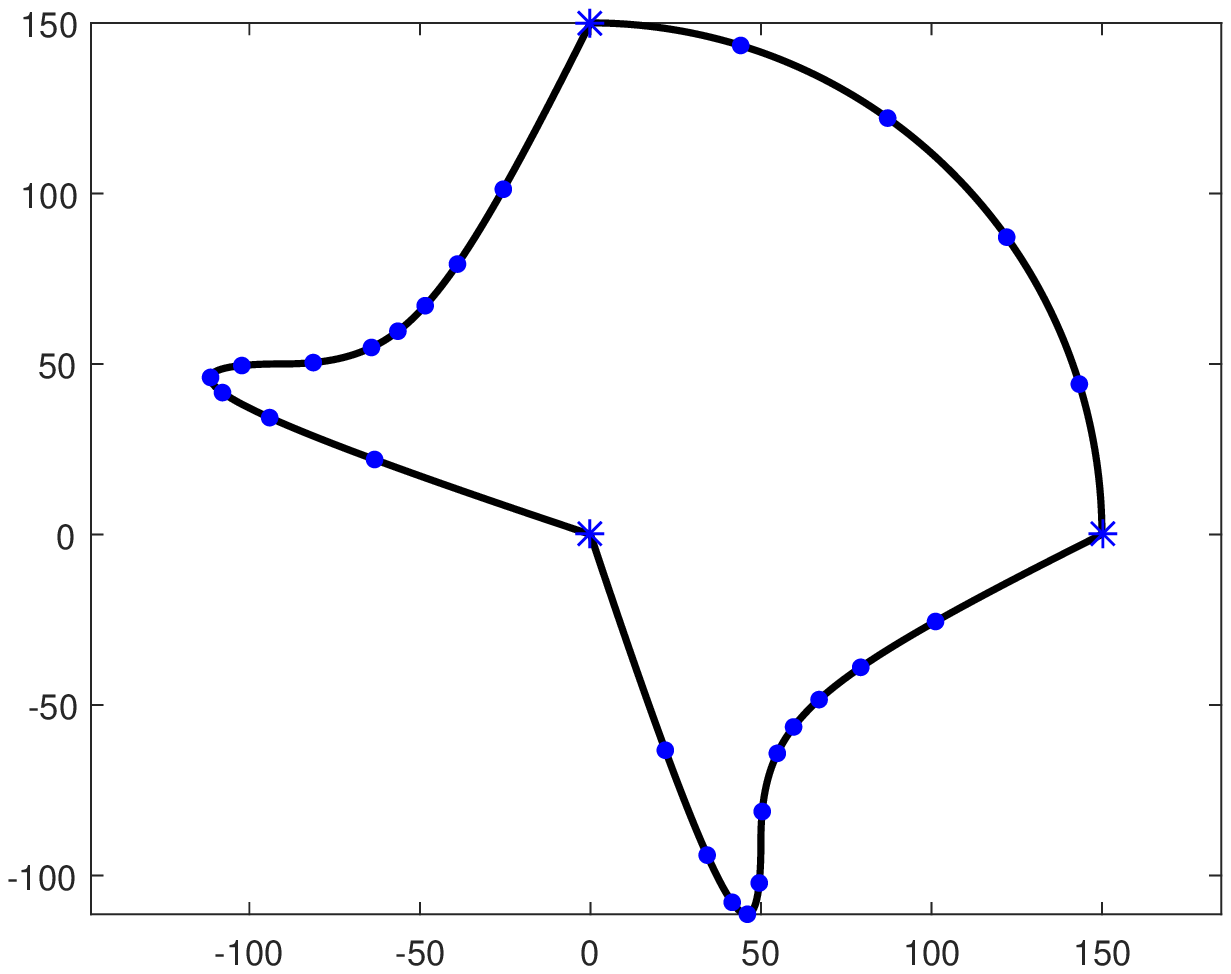}}
\subfigure[]{\includegraphics[scale=0.30]{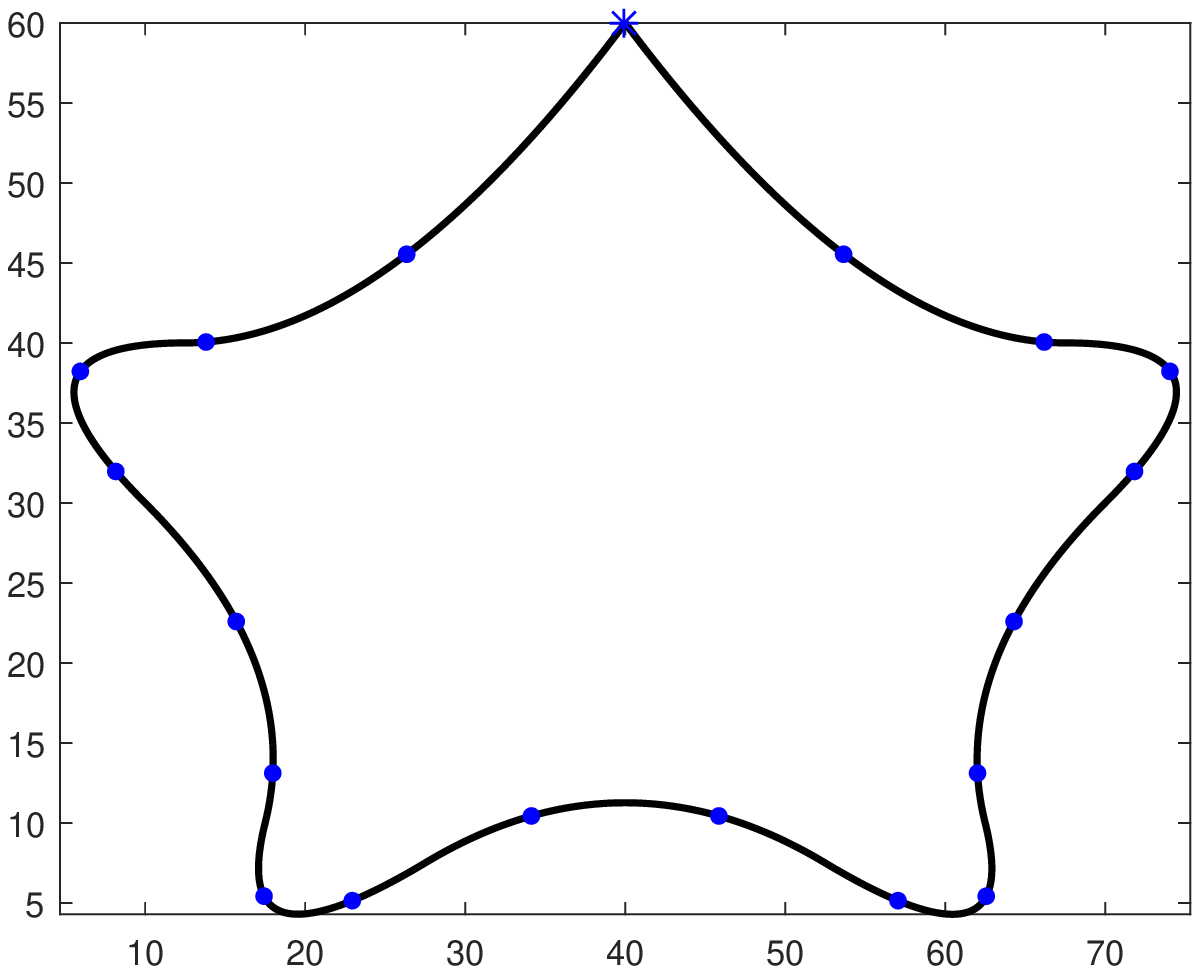}}
\subfigure[]{\includegraphics[scale=0.30, trim=-1.0cm -0.70cm 0cm 0cm,clip]{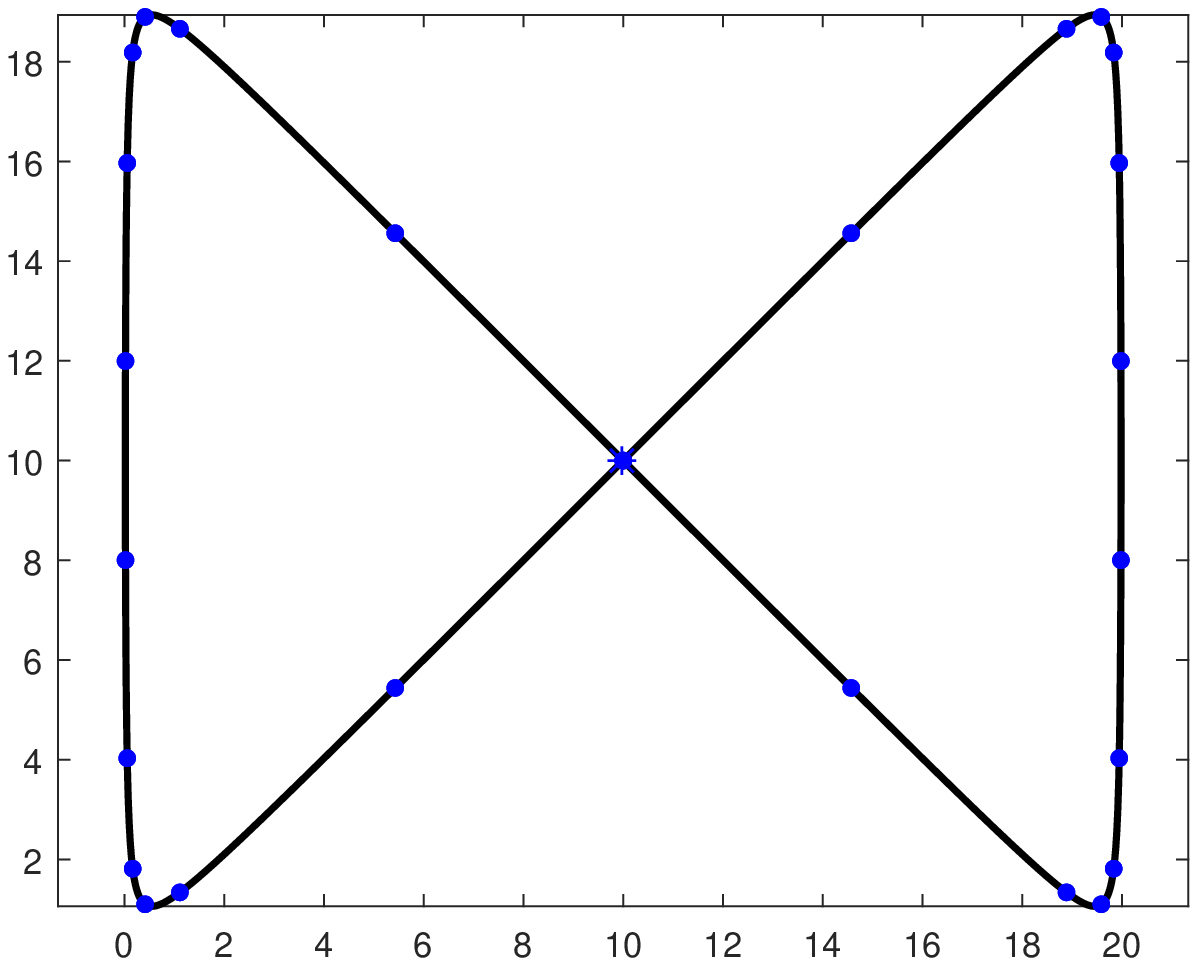}}\\
\subfigure[]{\includegraphics[scale=0.30]{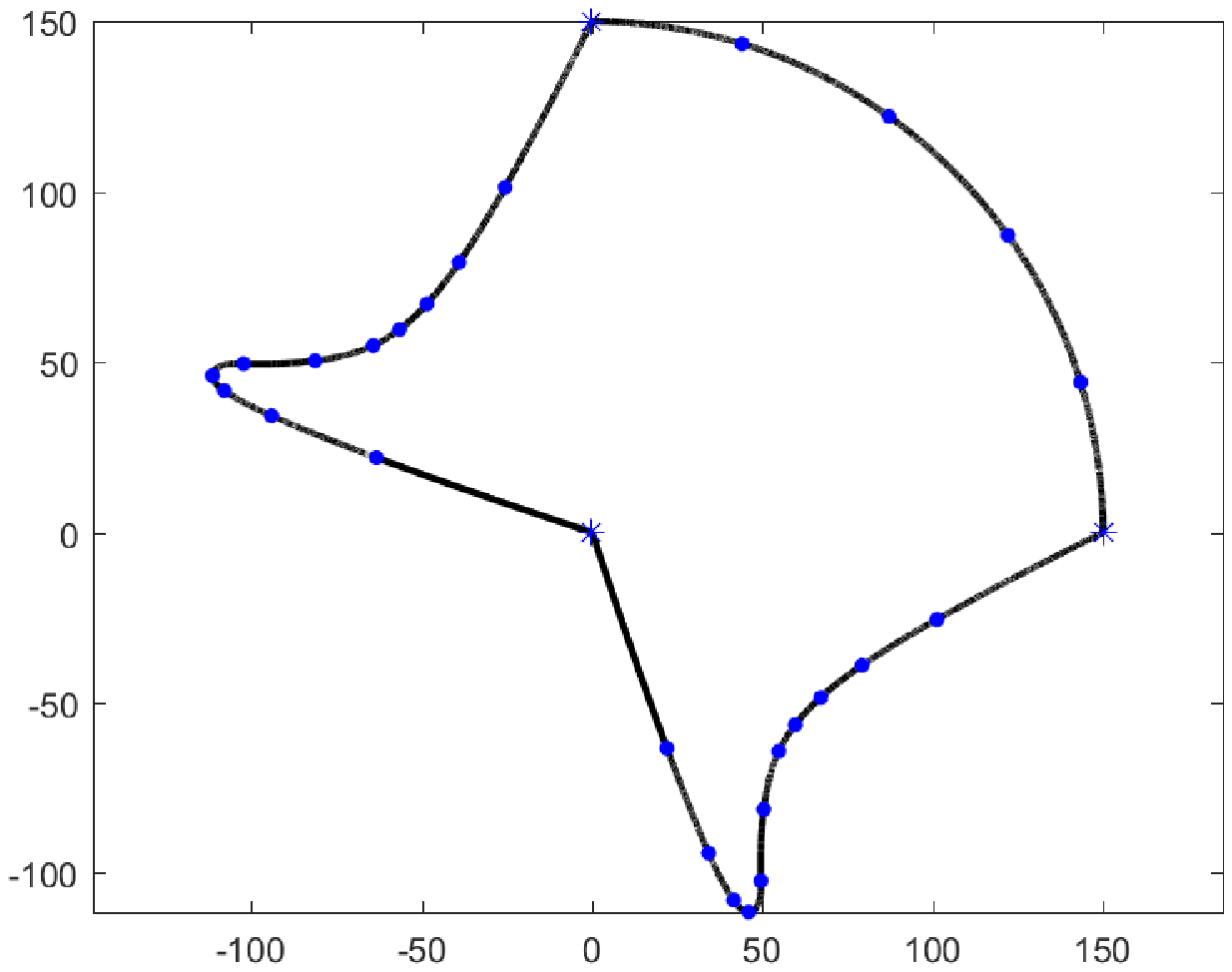}}
\subfigure[]{\includegraphics[scale=0.30]{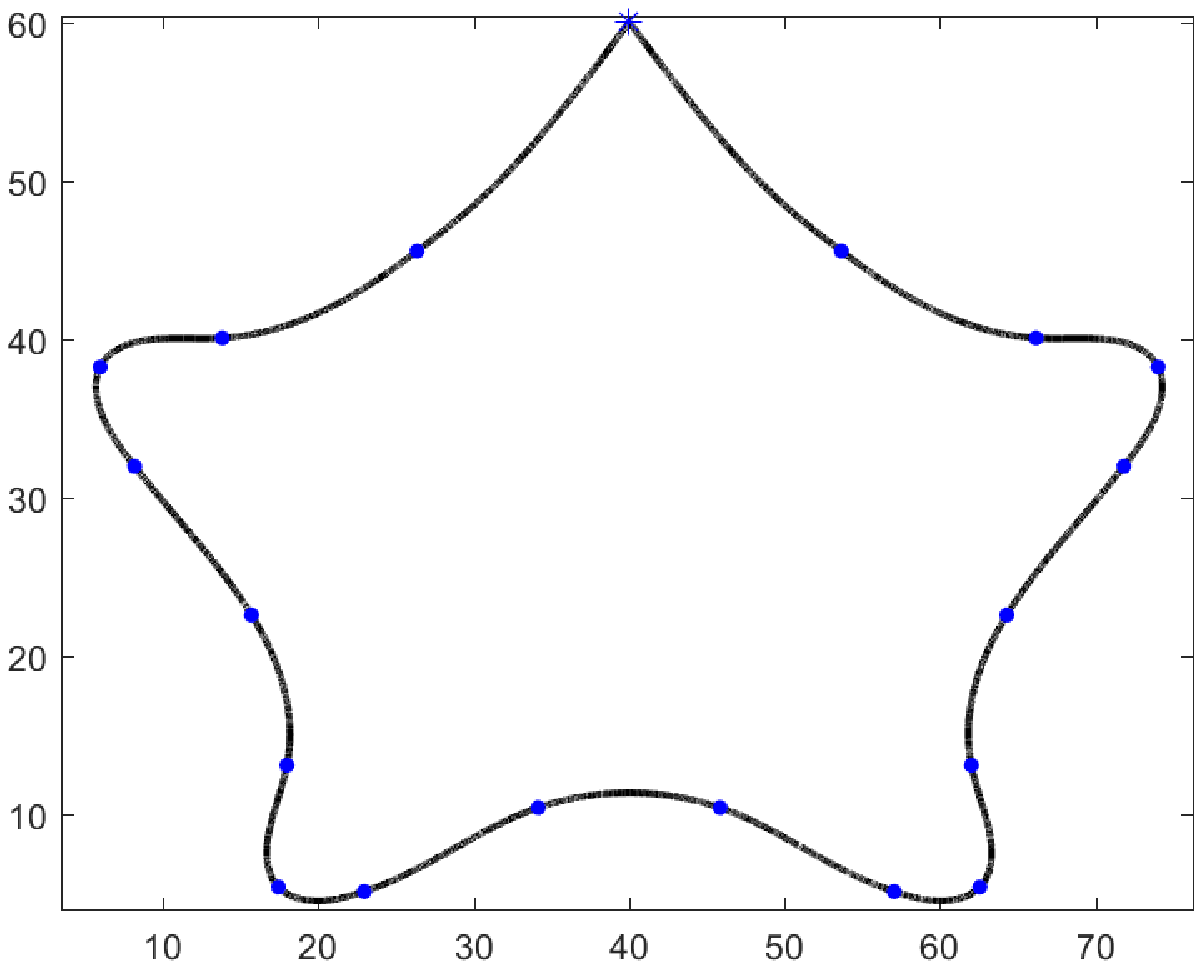}}
\subfigure[]{\includegraphics[scale=0.30]{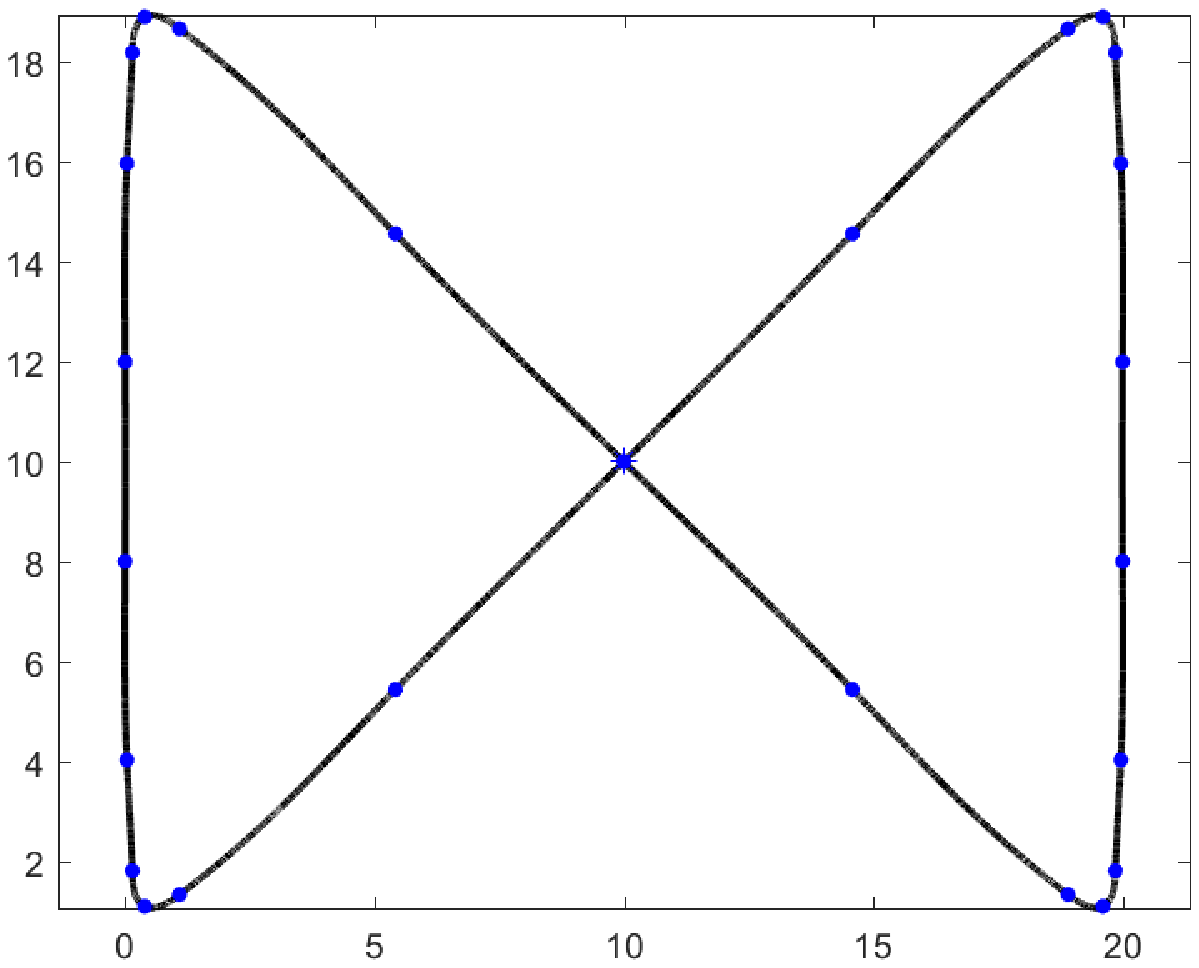}}
\caption{The three curved paths considered in the examples: the NURBS hat (a), the starfish (b), and the butterfly (c). Their PH spline approximations are shown in (d), (e), and (f), respectively.
The interpolation points are also shown (blue dots).}
\label{fig:NURBS}
\end{center}
\end{figure}

In this section three NURBS paths are considered to test the performance of the proposed feedrate scheduling algorithm. They are indicated as \textit{hat}, \textit{starfish}, and  \textit{butterfly} and were previously considered in \cite{lee11, xinhua16, annoni12}, respectively.  The three corresponding $C^2$ PH quintic spline approximations are shown in Figure \ref{fig:NURBS}  together with the original NURBS paths.
In particular, the hat curve is obtained with three PH spline segments interpolating $13, 6$, and $13$ points, respectively. For the other two curves instead, we used a single $C^2$ PH spline interpolating 20 (starfish) or 27 (butterfly) points on the original curve. These examples also confirm that PH curves are suitable for smooth approximation of any path.

The feedrate scheduler and the PH interpolator are implemented in Matlab and executed on a personal computer
with Intel(R) Core(TM) i5-4300U 1.90GHz CPU.  The feedrate scheduling parameters for the three examples are listed in Table \ref{tab:interp_param}, while the parameters $p_a$ and $p_j$ are always both set to $1/\sqrt{2}$ and $q_j$ to $1/2.$
The additional main parameters used in the experiments ---  the critical
curvature $\kappa_{cr}$, the number of special points, and the  simulation time --- are reported in Table \ref{tab:features} together with the results obtained by the scheduler with all
possible configurations of the relaxed and strict strategy. As expected, the strict strategy produces slightly lower motions.
Note also that the motion produced by this strategy is more segmented, since the number of crossing points is
higher than the number of critical points detected with the relaxed approach. When the simulation times reported in Table \ref{tab:features} are compared with those obtained with other methods, we may
	observe that our strategy can produce slightly slower motions. This is due to the fact that we always a priori ensures the desired acceleration control. 

\begin{table}[!h]
\centering
\caption{Parameter setting.}
\label{tab:interp_param}
\begin{tabular}{lrrrr}
\hline
\textbf{parameter}   & \textbf{symbol}   & \textbf{hat}   &
\textbf{starfish}   & \textbf{butterfly} \\ \hline 
commanded
feedrate   (mm/s)  		 & $V_m$      & $250$ & $200$ & $100$ \\
maximum acceleration (mm/s$^2$) & $A_m$      & $800$  & $3000$ & $1000$ \\ 
maximum jerk (mm/s$^3$)       & $J_m$      & $26400$ & $60000$ & $50000$ \\ 
maximum chord error  (mm) & $D$  	  & $0.0010$ & $0.0010$ & $0.0002$\\
sampling time (s) & $\Delta t$ & $0.002$  & $0.001$ & $0.0005$\\\hline
\end{tabular}
\end{table}

\begin{table}[!h]
\centering
\caption{Results of the scheduler obtained with the relaxed ($R_0, R_1, R_2$) and strict strategies ($S_0, S_1, S_2$).}
\label{tab:features}
\begin{tabular}{lrrr}
\hline
  & \textbf{hat}   & \textbf{starfish}   &
\textbf{butterfly} \\ 
\hline
$\kappa_{cr}$ (mm$^{-1}$)     & $9.9561 \cdot 10^{-3}$   & $5.3033 \cdot 10^{-2}$ 		& $9.051 \cdot 10^{-3}$\\ 
\hline
\multicolumn{4}{c}{\emph{relaxed strategy}}\\
\hline
\texttt{\#} critical points  & $4$      & $10$		& 4 \\
${R_0}$ simulation time (s)      & $7.016$ s  & $1.396$	& $1.854$  \\ 
$R_1$ simulation time (s)    & $7.045$  & $1.453$	& $1.899$ \\ 
$R_2$ simulation time (s)     & $7.159$ & $1.903$ & $2.474$ \\ 
\hline
\multicolumn{4}{c}{\emph{strict strategy}} \\
\hline
\texttt{\#} crossing points  & $8$      & $20$		& $12$ \\
{$S_0$} simulation time (s)      & $7.820$  & $1.749$ & $2.237$\\
{$S_1$} simulation time (s)      & $7.884$  & $2.947$	 & $2.405$  \\
{$S_2$} simulation time (s)      & $9.141$  & $3.230$ & $4.810$ \\
\hline
\end{tabular}
\end{table}

All the following figures corresponds to the results obtained with the relaxed formulation of the scheduler, since there was no need to consider the strict formulation to satisfy the prescribed bounds in these cases. In particular we show the results obtained with $R_0$ for the \textit{hat curve}, $R_1$ for the \textit{starfish} example, and $R_1$
for the \textit{butterfly} path.

Unlike \cite{lee11}, we are able to prescribe a $C^2$ feedrate function. Note that the feedrate defined in \cite{xinhua16} is also $C^2$ but with a less smooth profile which can produce a major vibration during the machining. In the presented results we
obtain acceleration continuous motions since we apply our feedrate scheduler to $C^2$ PH splines whose curvature first derivative has generally a finite jump at the spline knots. Jerk continuity could be obtained  by applying the scheduler to curves with continuous curvature derivative.

\begin{figure}[!h]
\begin{center}
\includegraphics[scale=0.33]{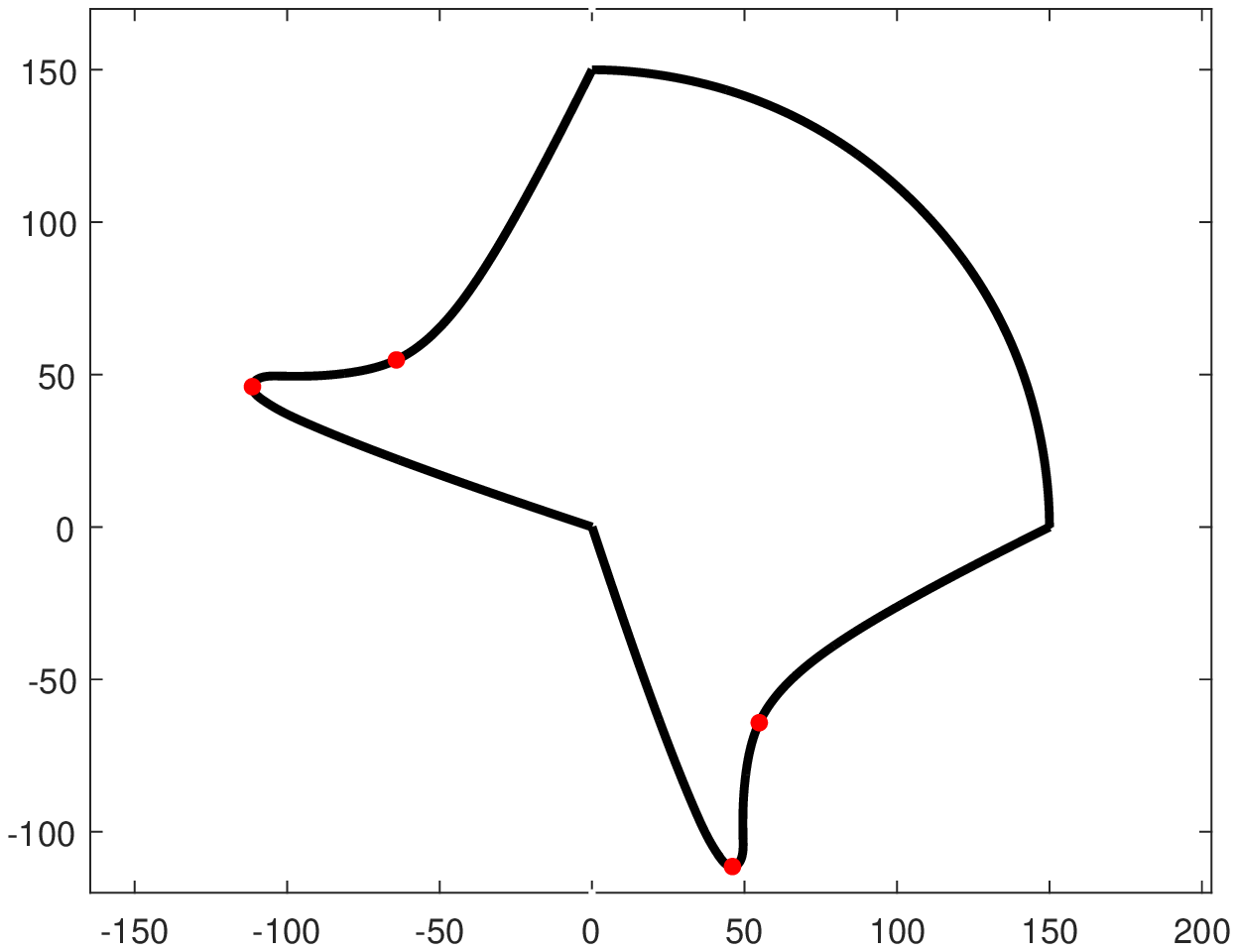}
\includegraphics[scale=0.33,trim=0cm 0.75cm 0cm 0cm,clip]{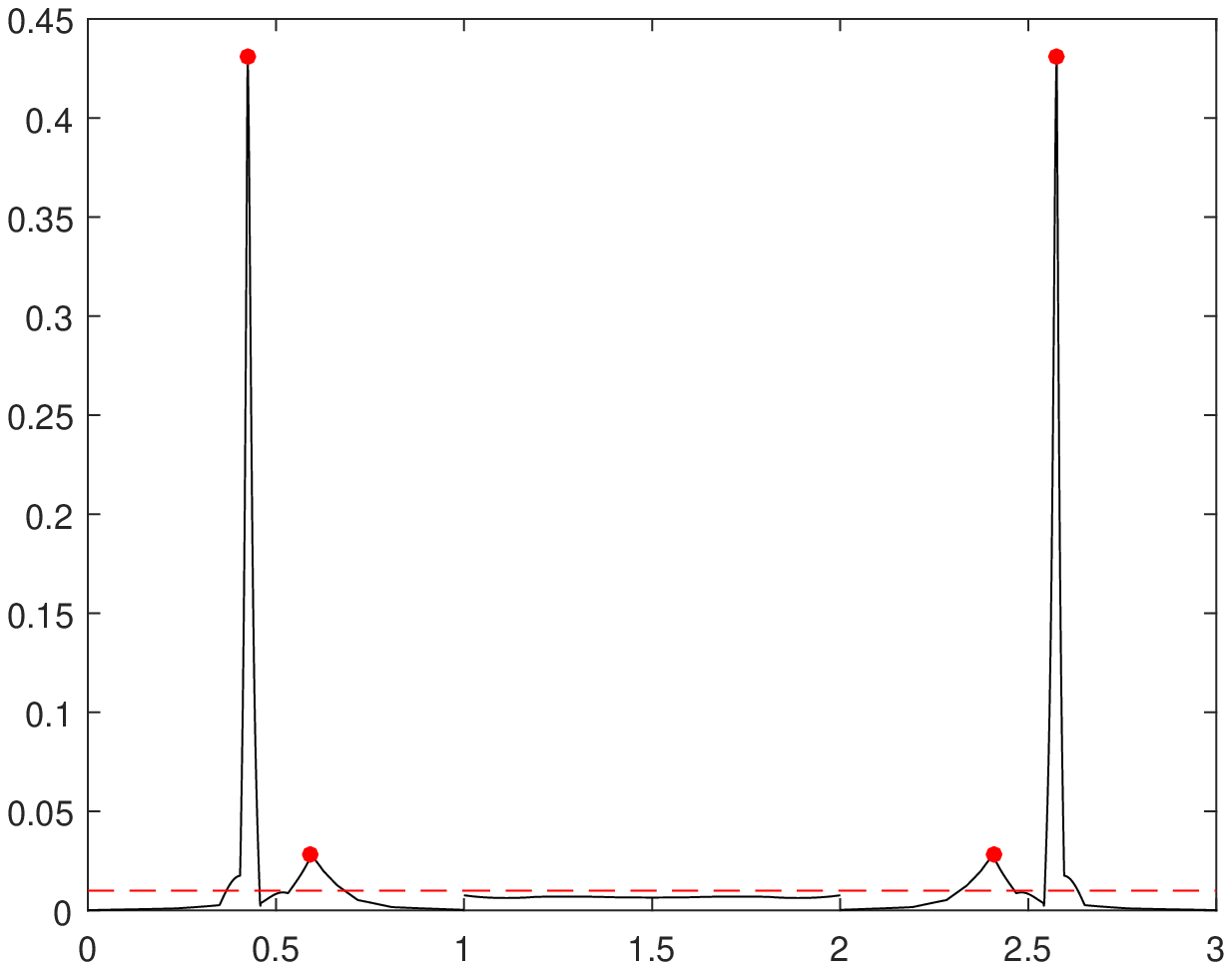}
\includegraphics[scale=0.33,trim=1.1cm 0.69cm -1cm 0cm,clip]{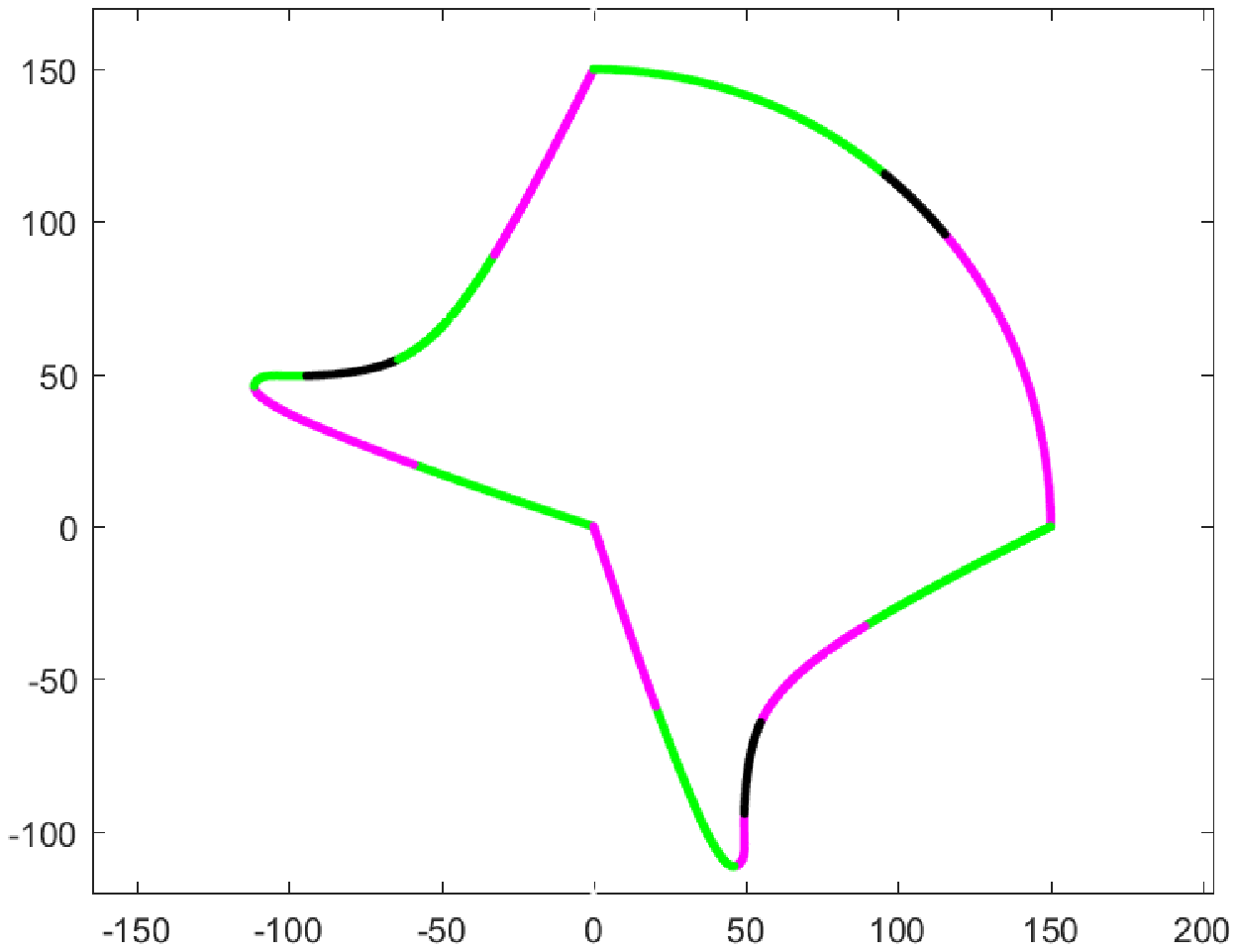}
\caption{Hat curve example: PH spline approximation (left), absolute curvature plot and critical curvature value (center, black and dashed red lines), and reference points (right). The critical points on the curve and the corresponding peaks on the
curvature plot are also shown (red dots).}
\label{fig:hat}
\end{center}
\end{figure}

\begin{figure}[!h]
\begin{center}
\subfigure[]{\includegraphics[scale=0.40,trim=0cm 0cm -1cm 0cm,clip]{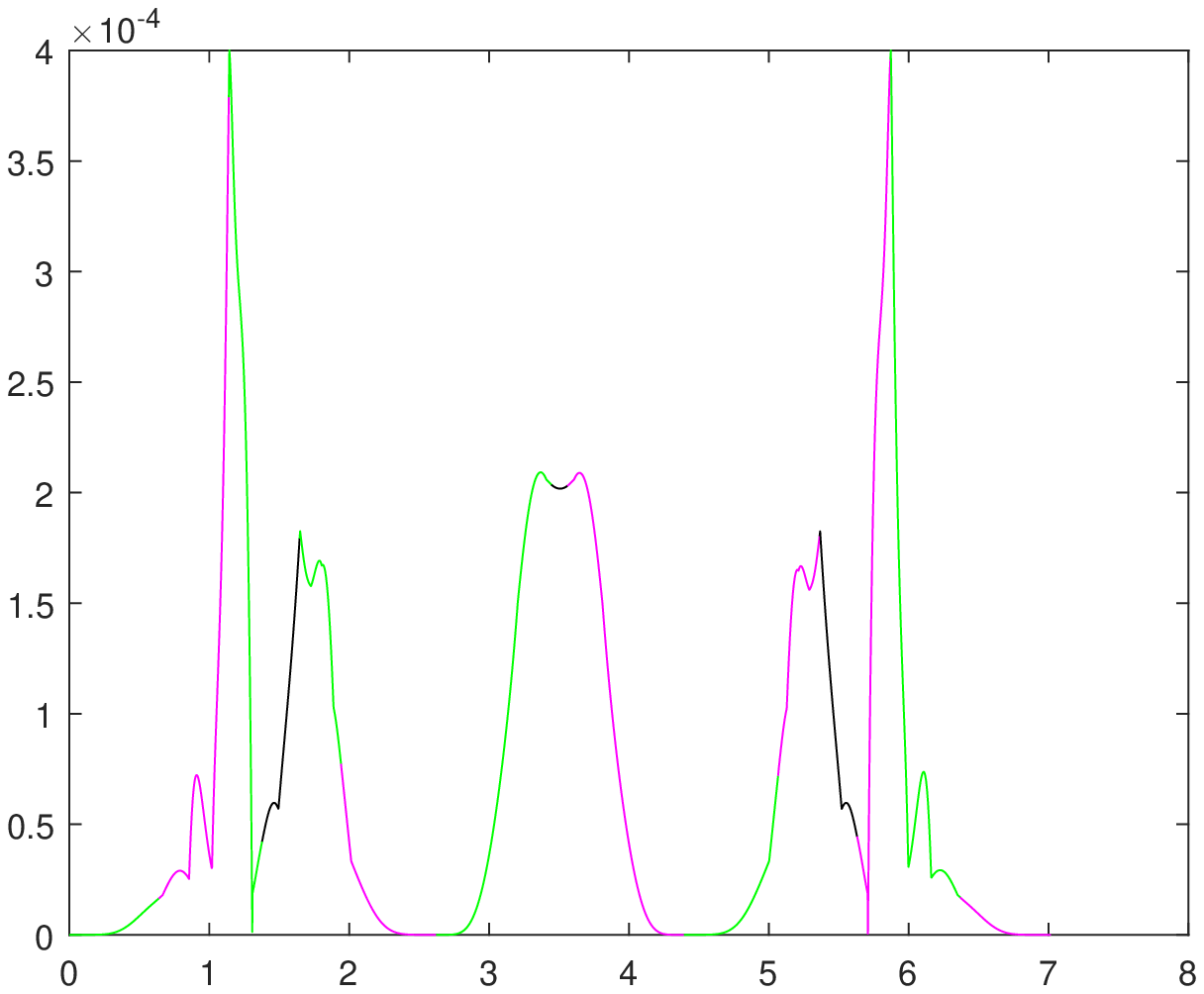}}\ \ \
\subfigure[]{\includegraphics[scale=0.40]{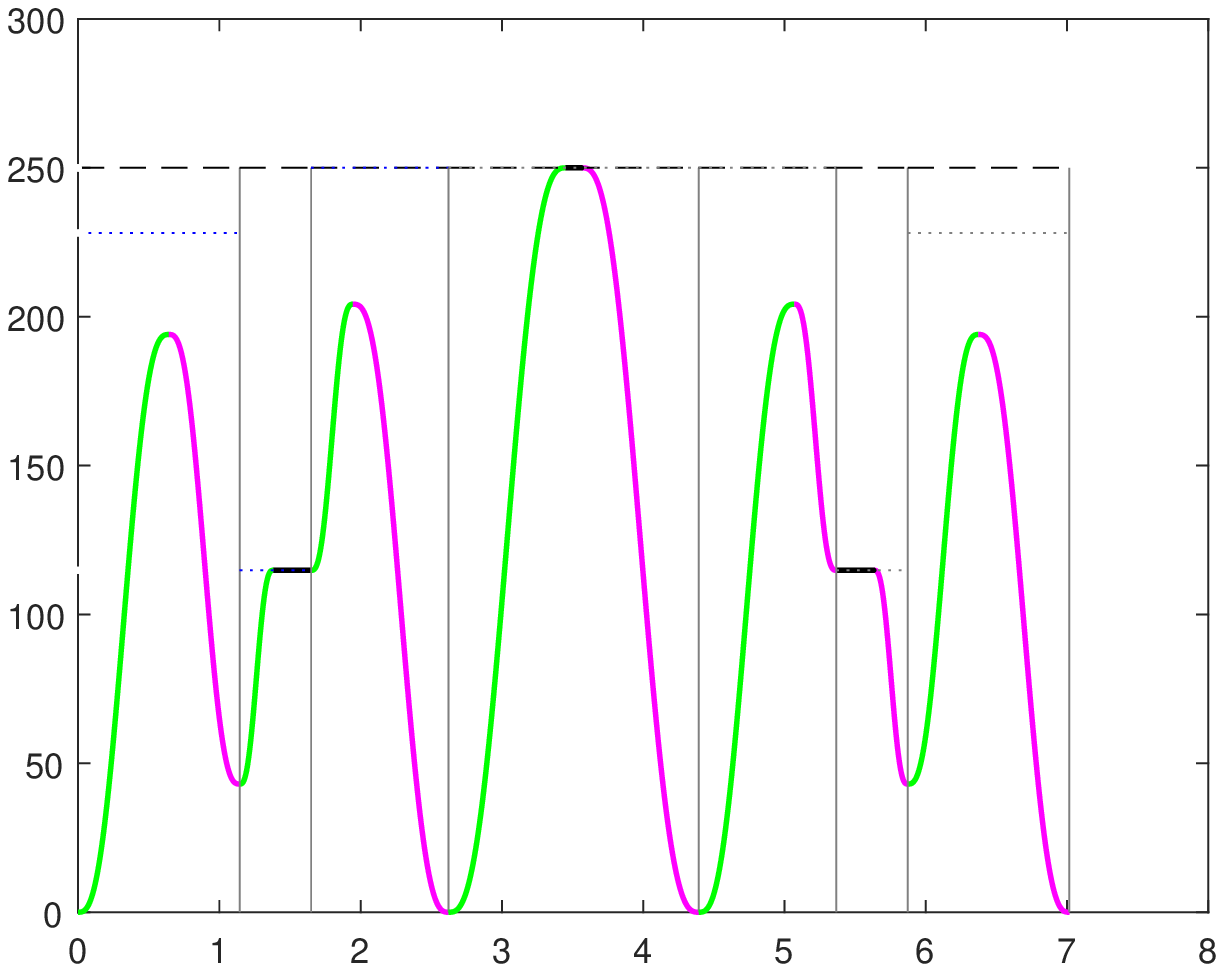}}\\
\subfigure[]{\includegraphics[scale=0.40,trim=0cm 0cm -1cm 0cm,clip]{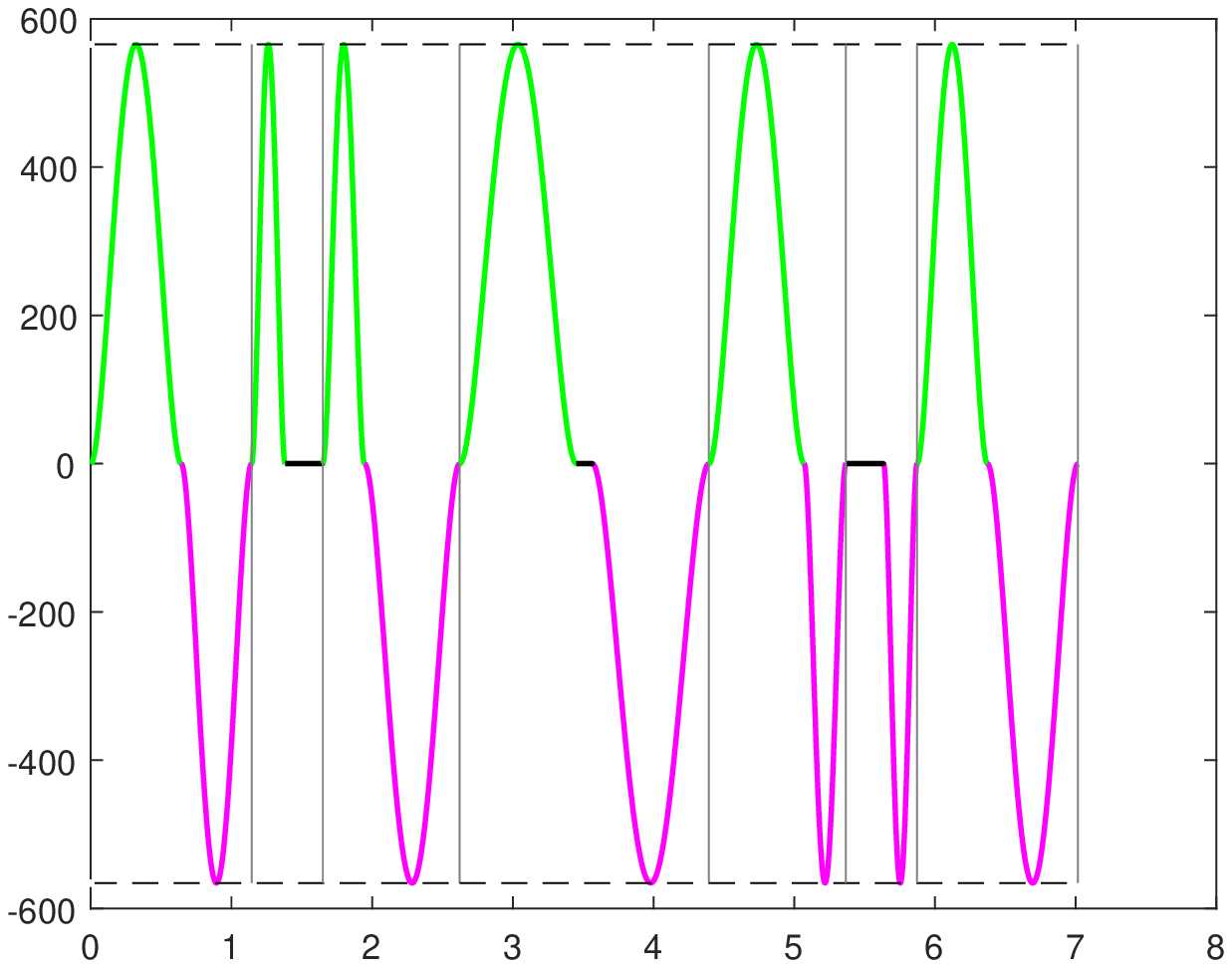}}\ \ \
\subfigure[]{\includegraphics[scale=0.40]{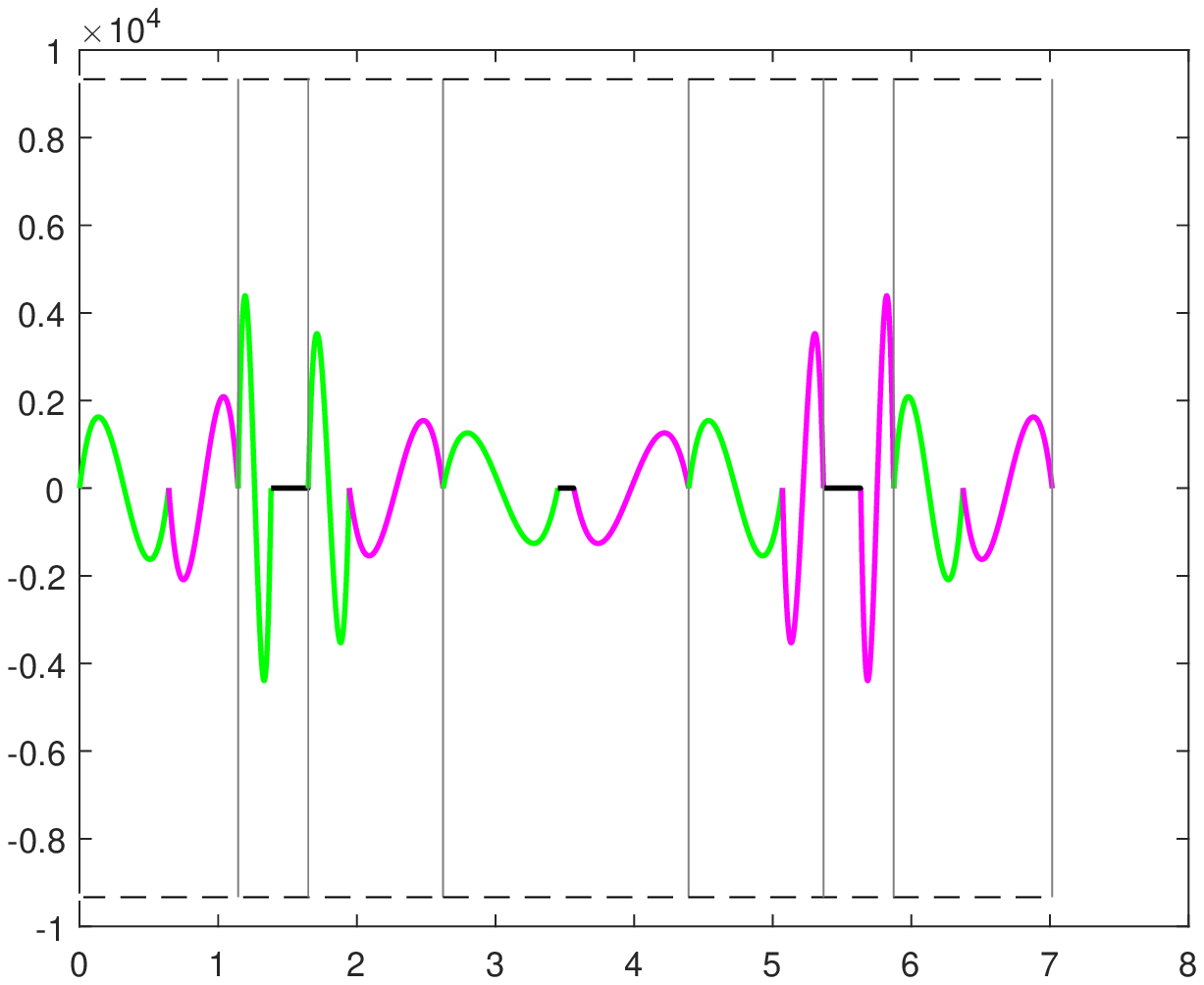}}
\caption{Feedrate scheduling results for the hat curve obtained using the $R_0$ configuration. The chord error (a) and the feedrate profile (b) are shown together with the first (c) and second (d) derivatives of the feedrate function.}
\label{fig:hat_feedrate}
\end{center}
\end{figure}

\begin{figure}[h]
\begin{center}
\subfigure[]{\includegraphics[scale=0.40,trim=0cm 0cm -1cm 0cm,clip]{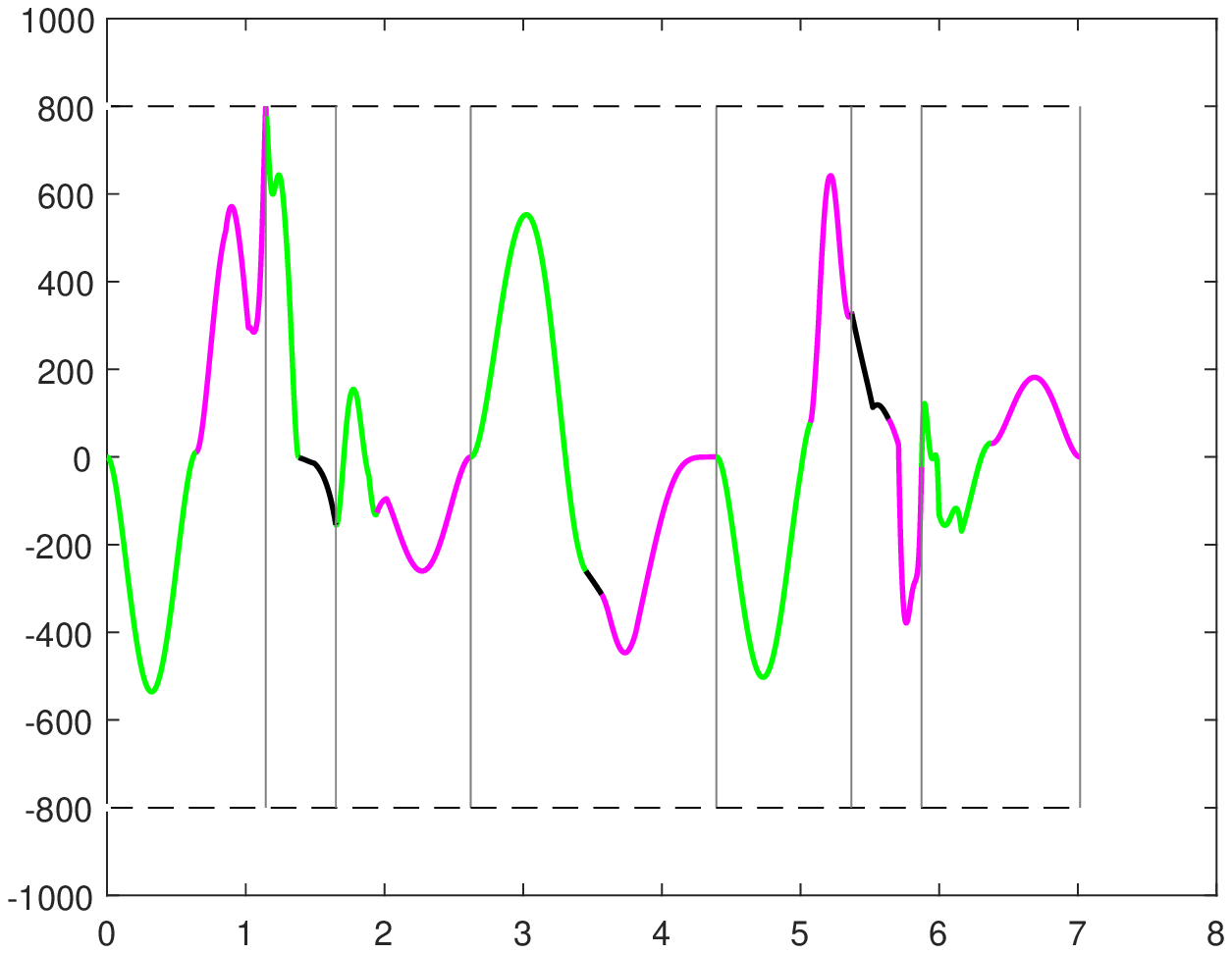}}\ \ \
\subfigure[]{\includegraphics[scale=0.40]{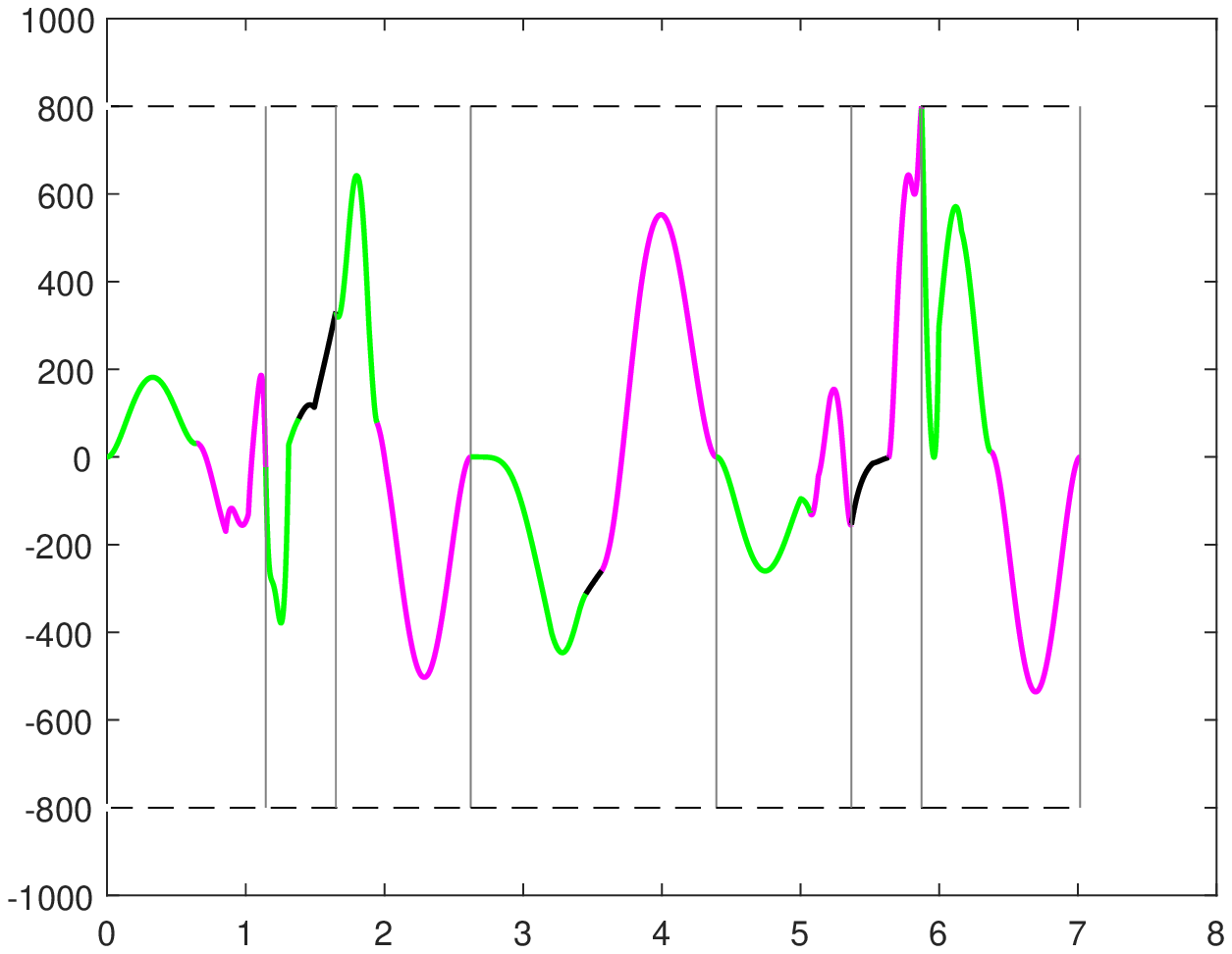}}\\
\subfigure[]{\includegraphics[scale=0.40,trim=0cm 0cm 1.3cm 0cm,clip]{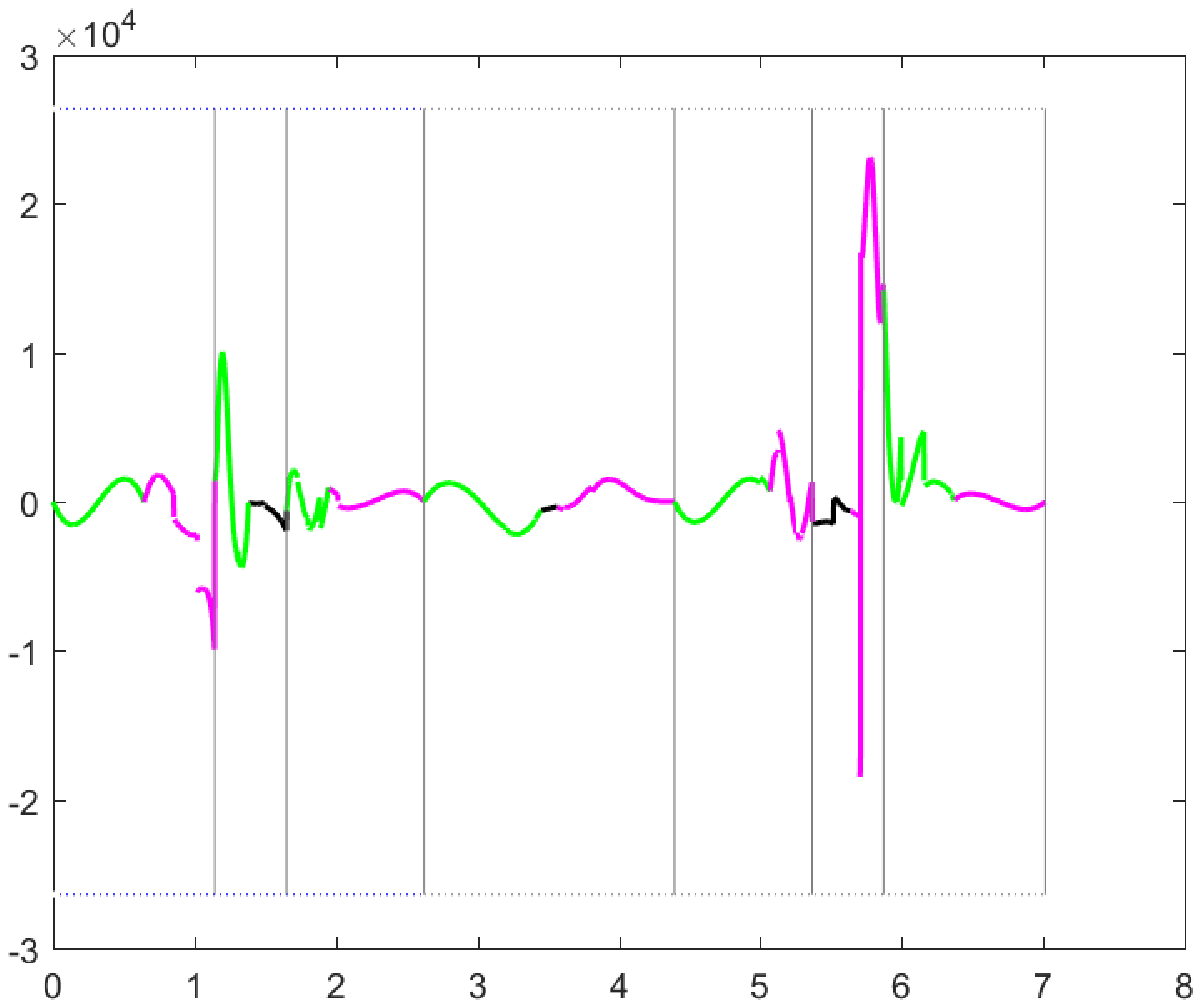}}\ \ \
\subfigure[]{\includegraphics[scale=0.40]{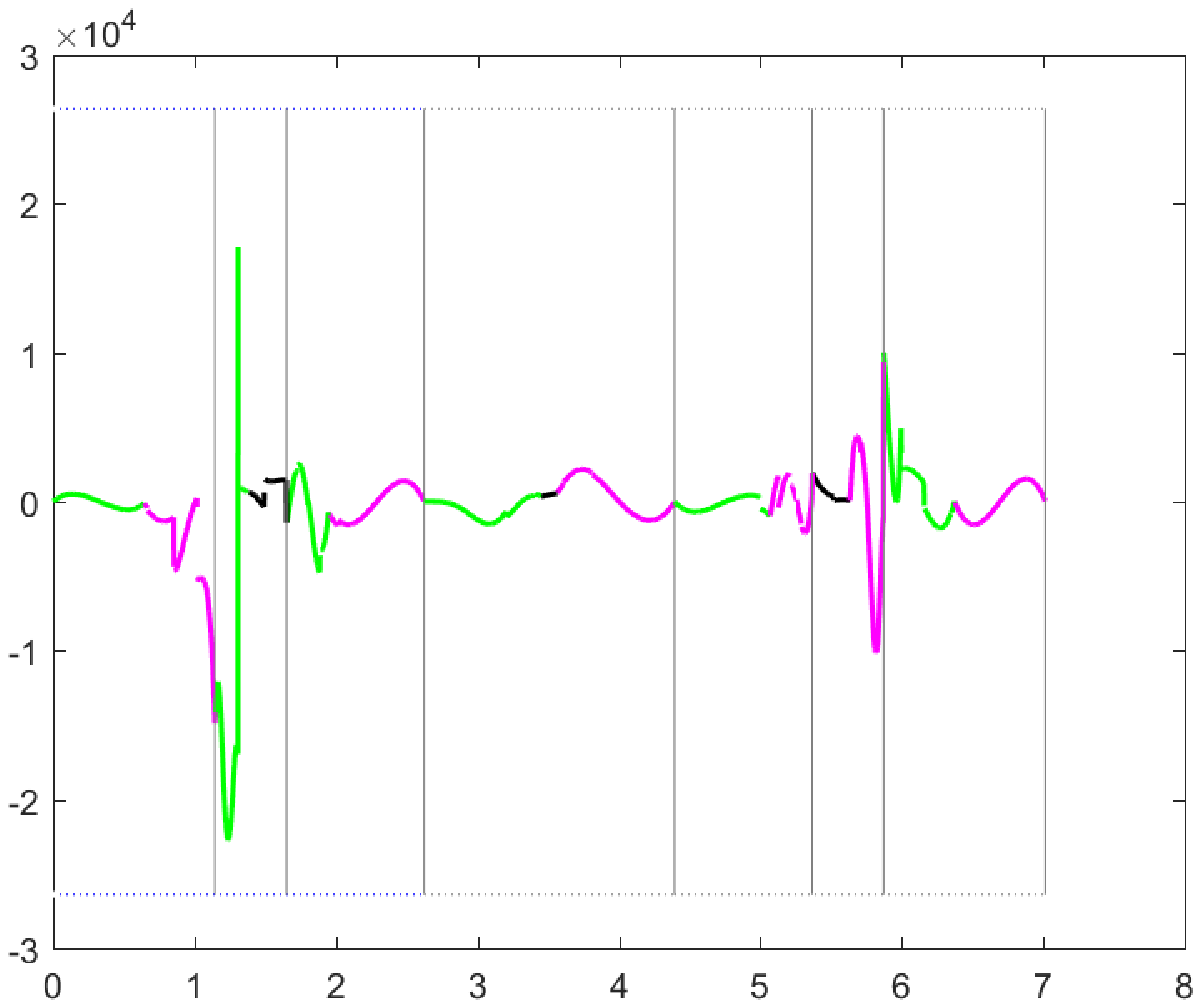}} 
\caption{Cartesian components of acceleration (a-b) and jerk (c-d)  obtained for the hat curve with the $R_0$ configuration.}
\label{fig:hat_cartesian}
\end{center}
\end{figure}

\begin{figure}[!h]
\begin{center}
\includegraphics[scale=0.33]{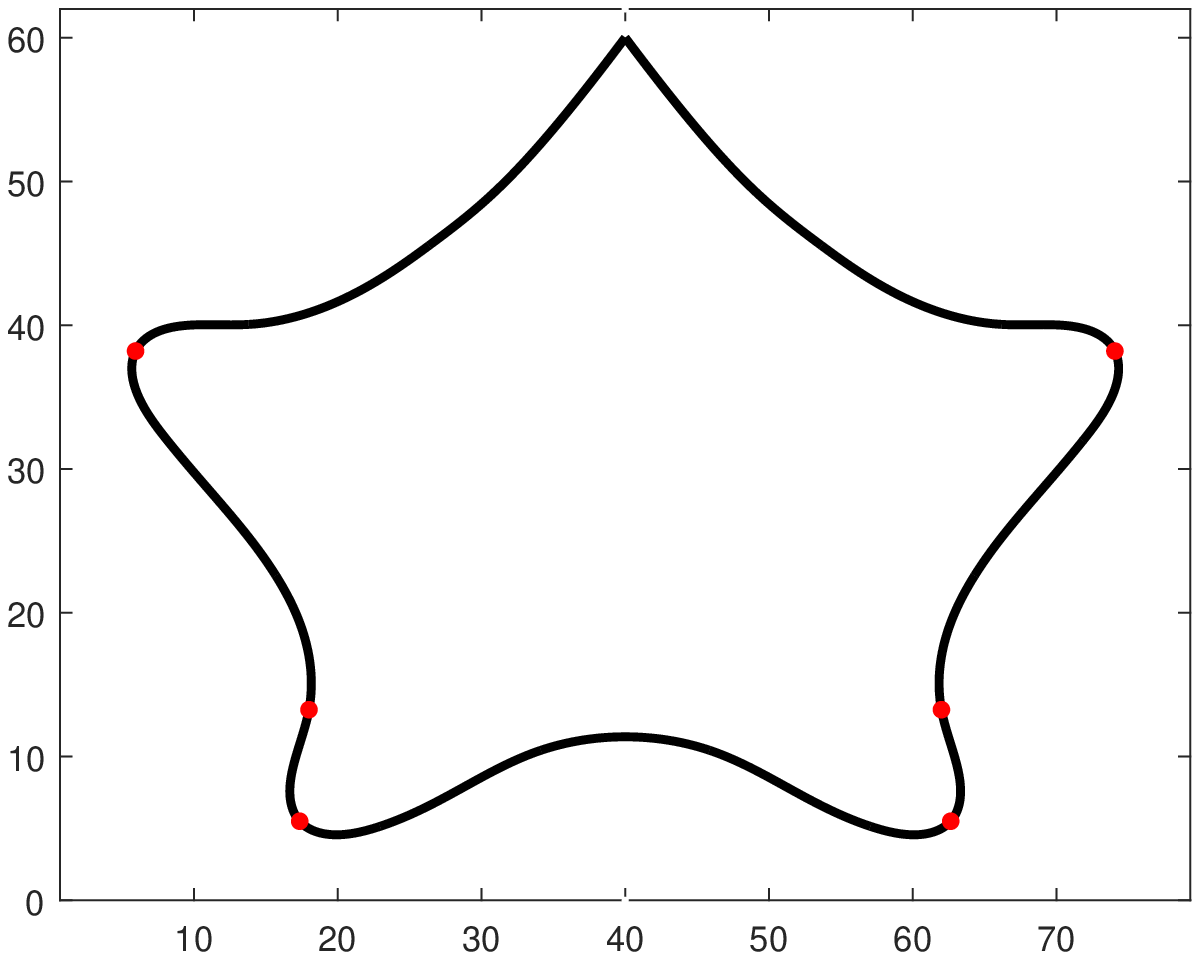}
\includegraphics[scale=0.33,trim=0cm 0.0cm 0cm 0cm,clip]{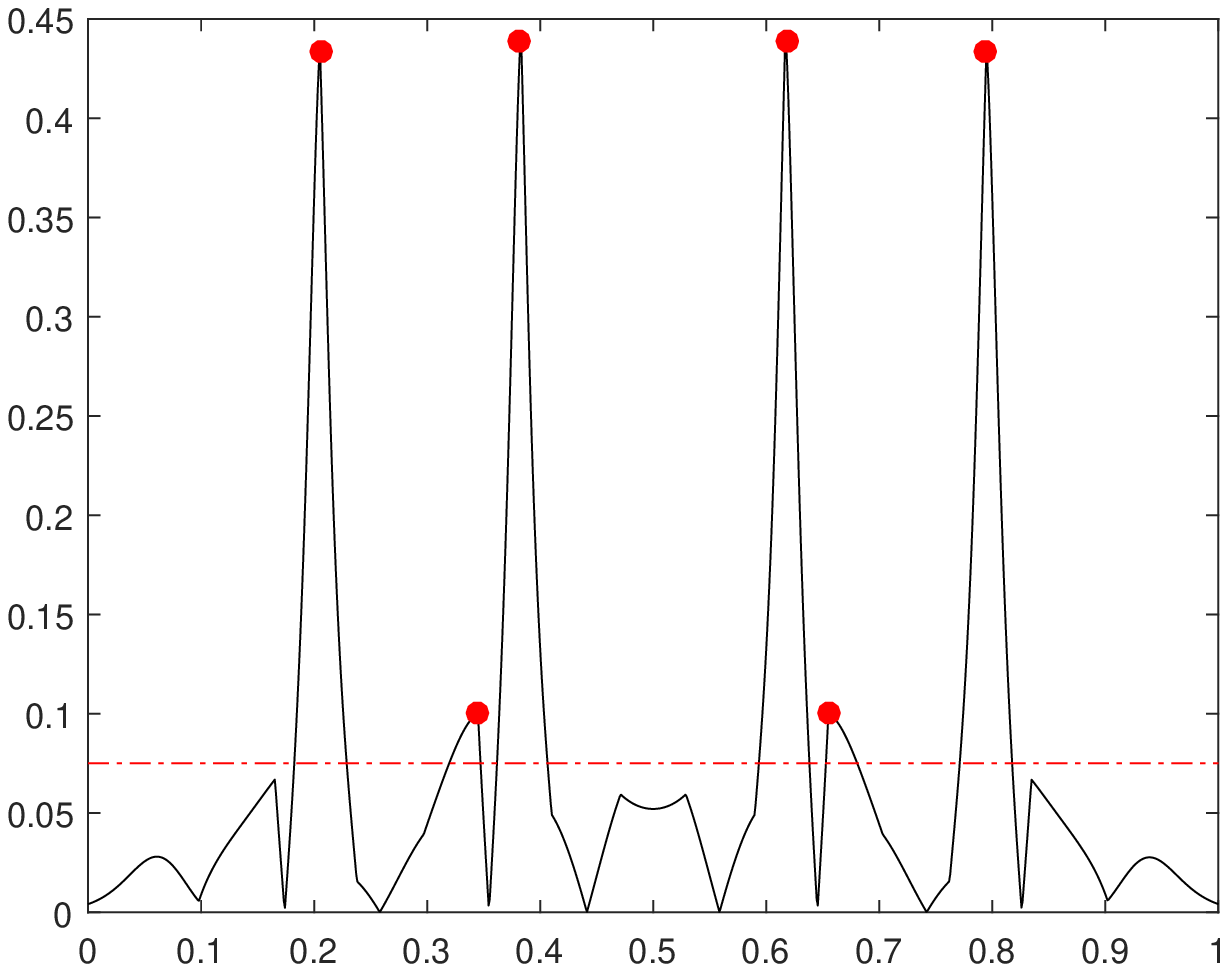}
\includegraphics[scale=0.33,trim=1.1cm 0.69cm -1cm 0cm,clip]{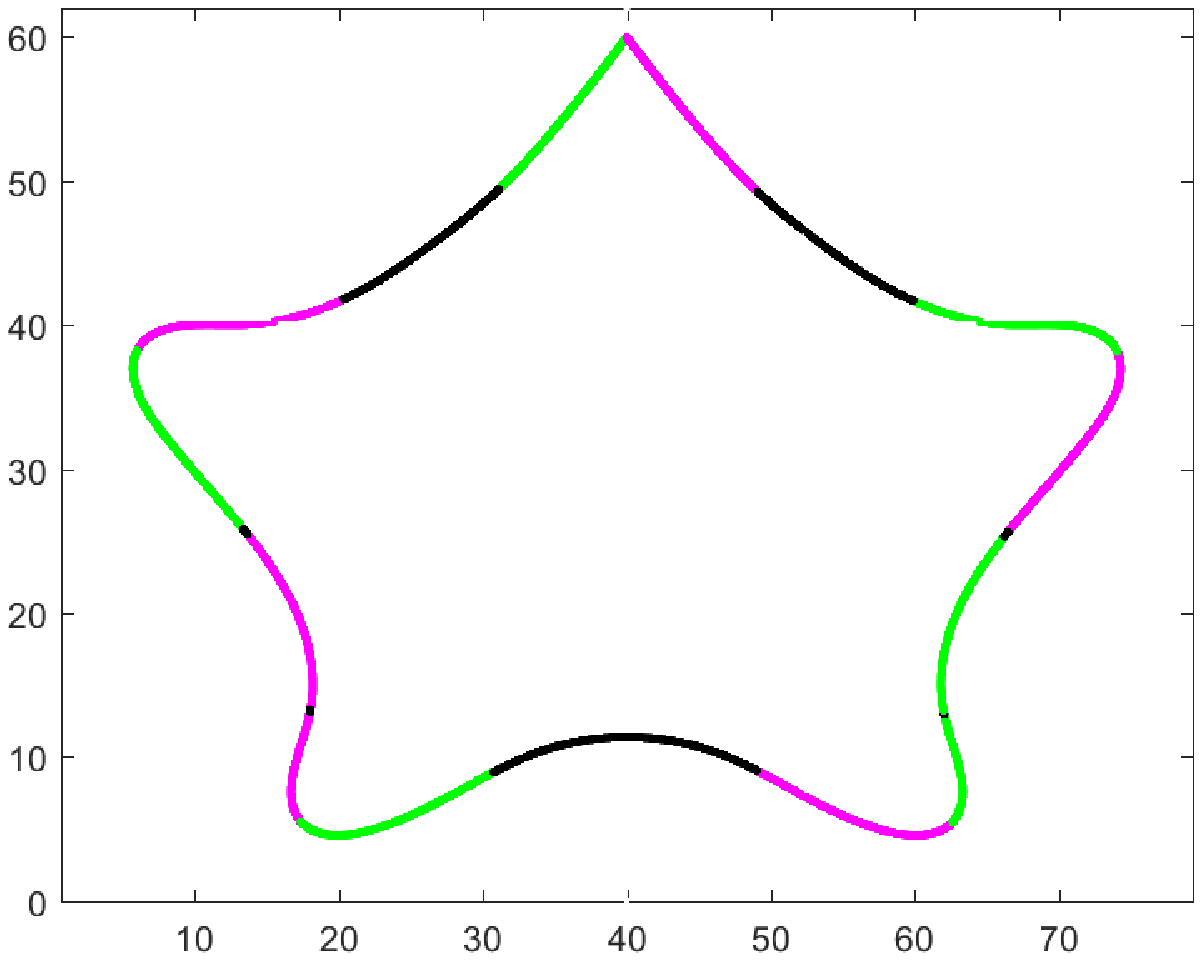}
\caption{Starfish curve example: PH spline approximation (left), absolute curvature plot and critical curvature value (center, black and dashed red lines), and reference points (right). The critical points on the curve and the corresponding peaks on the
curvature plot are also shown (red dots).}
\label{fig:star}
\end{center}
\end{figure}

\begin{figure}[h]
\begin{center}
\subfigure[]{\includegraphics[scale=0.40,trim=0cm 0cm -1cm 0cm,clip]{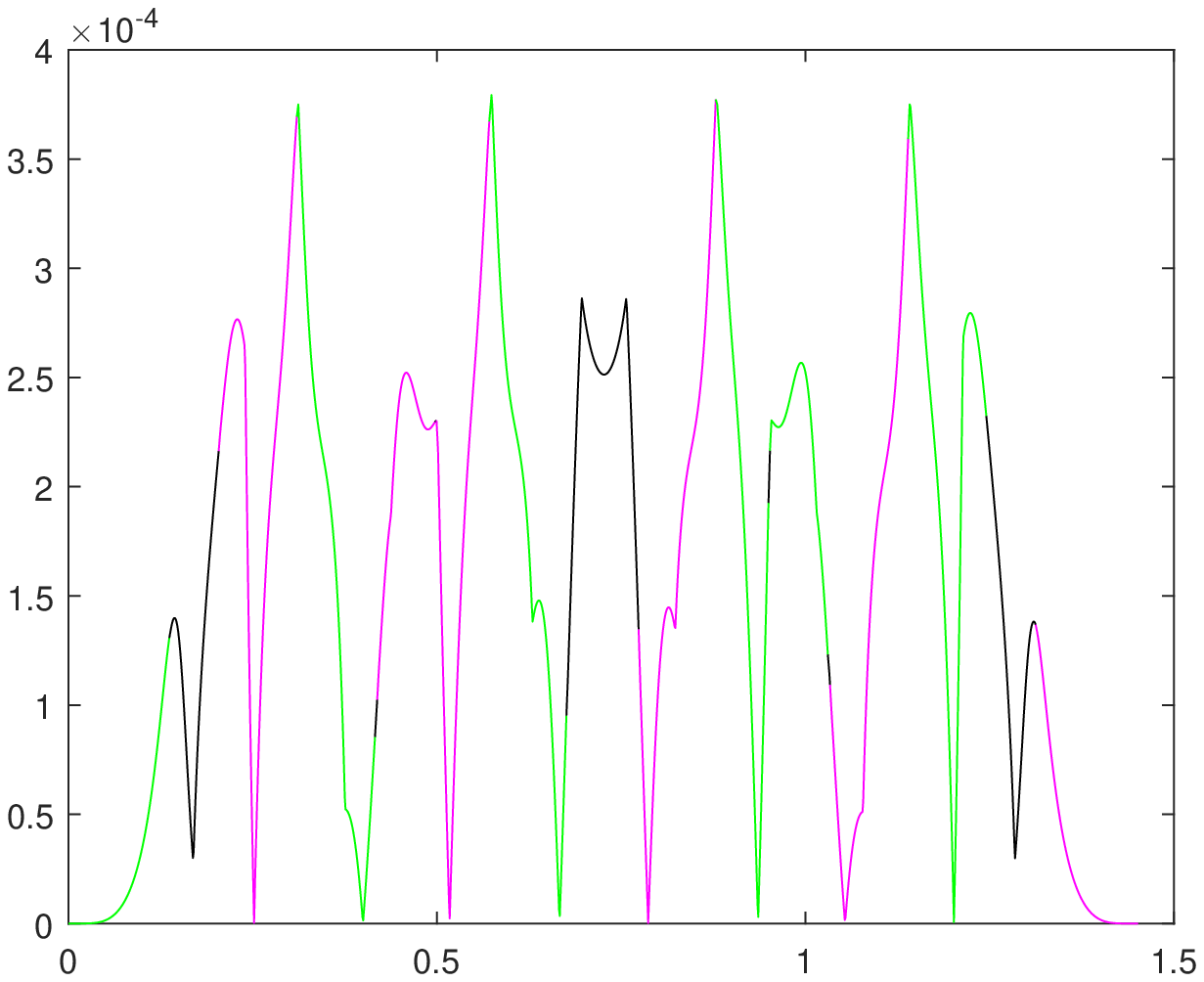}}\ \ \
\subfigure[]{\includegraphics[scale=0.40]{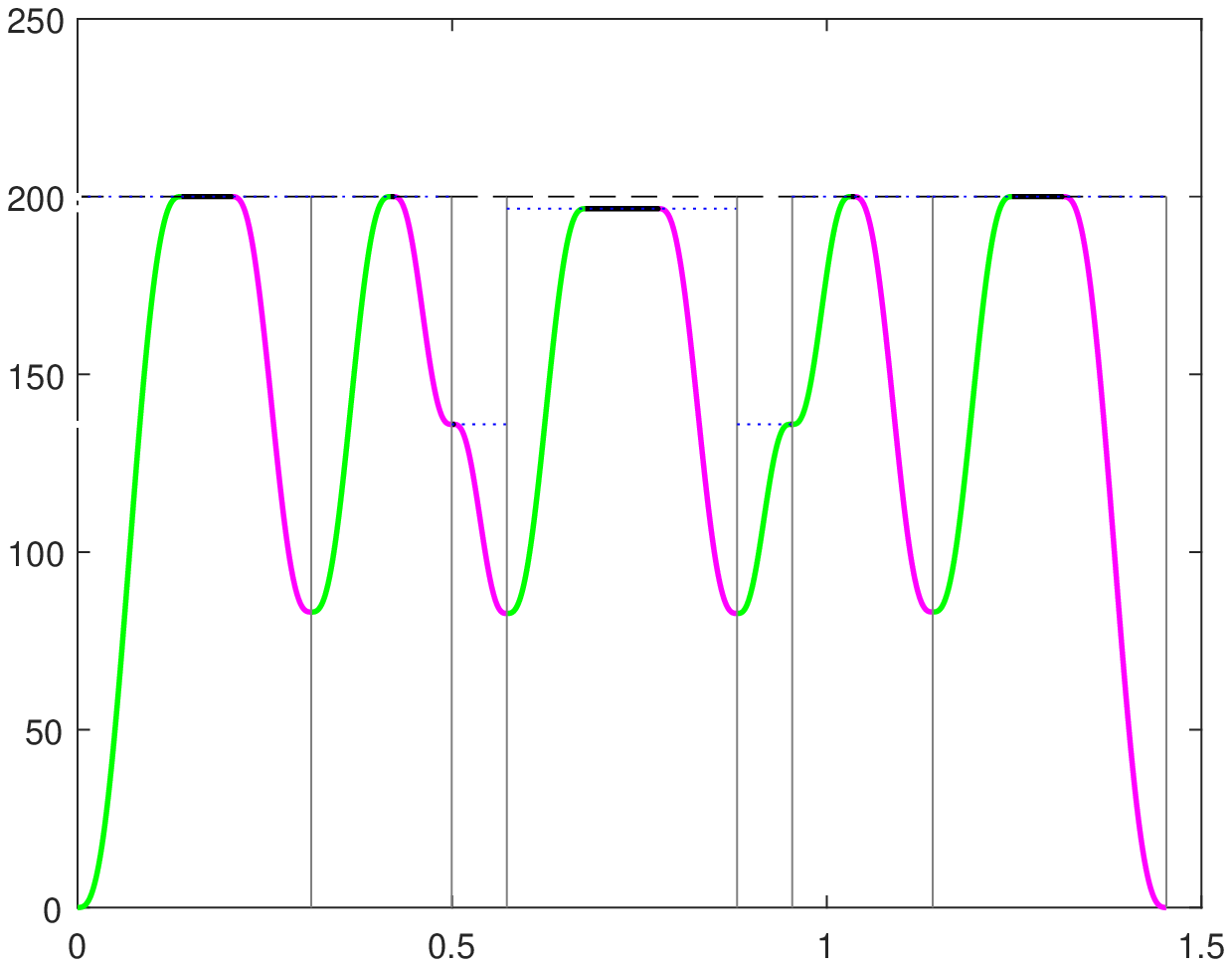}}\\
\subfigure[]{\includegraphics[scale=0.40,trim=0cm 0cm -1cm 0cm,clip]{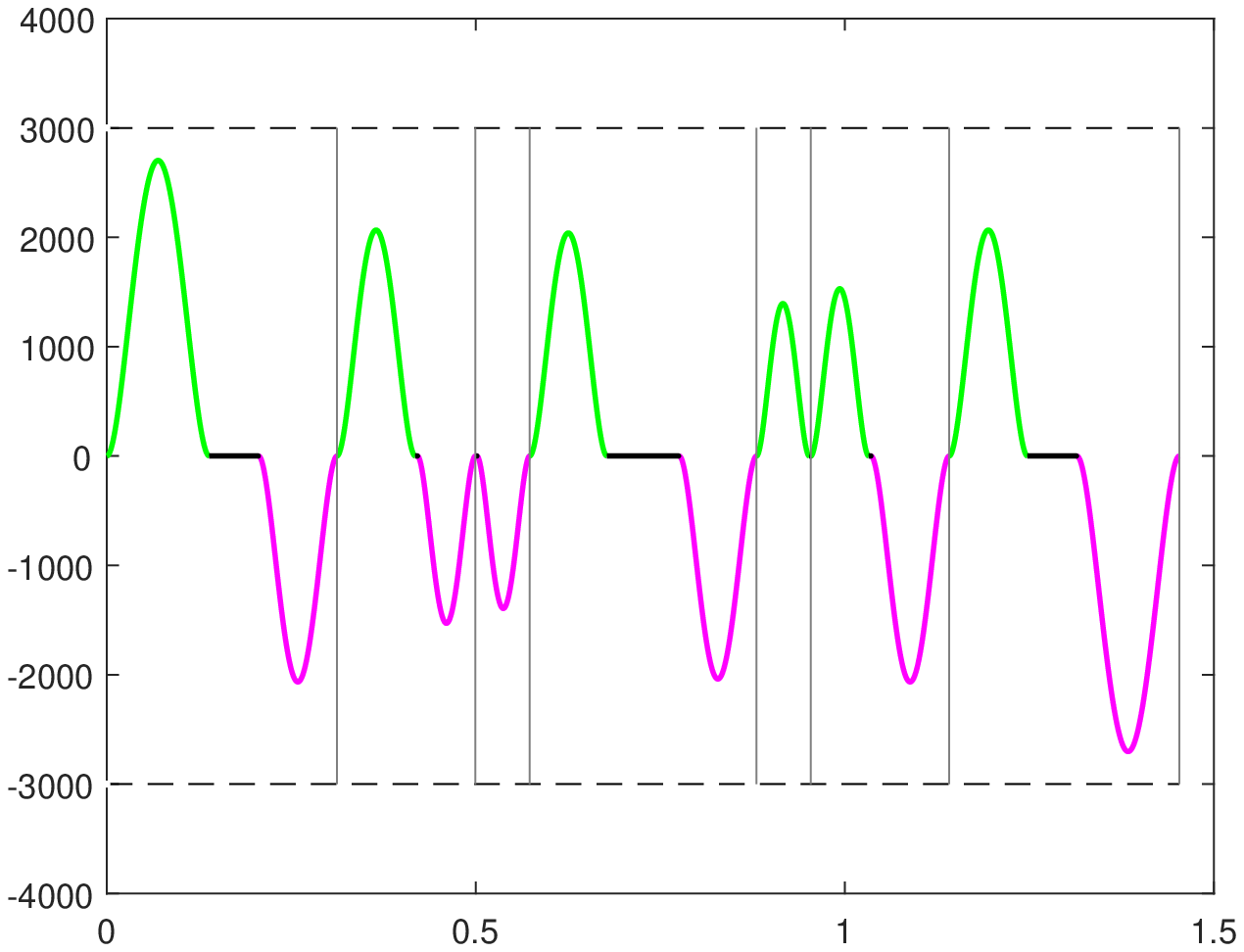}}\ \ \
\subfigure[]{\includegraphics[scale=0.40]{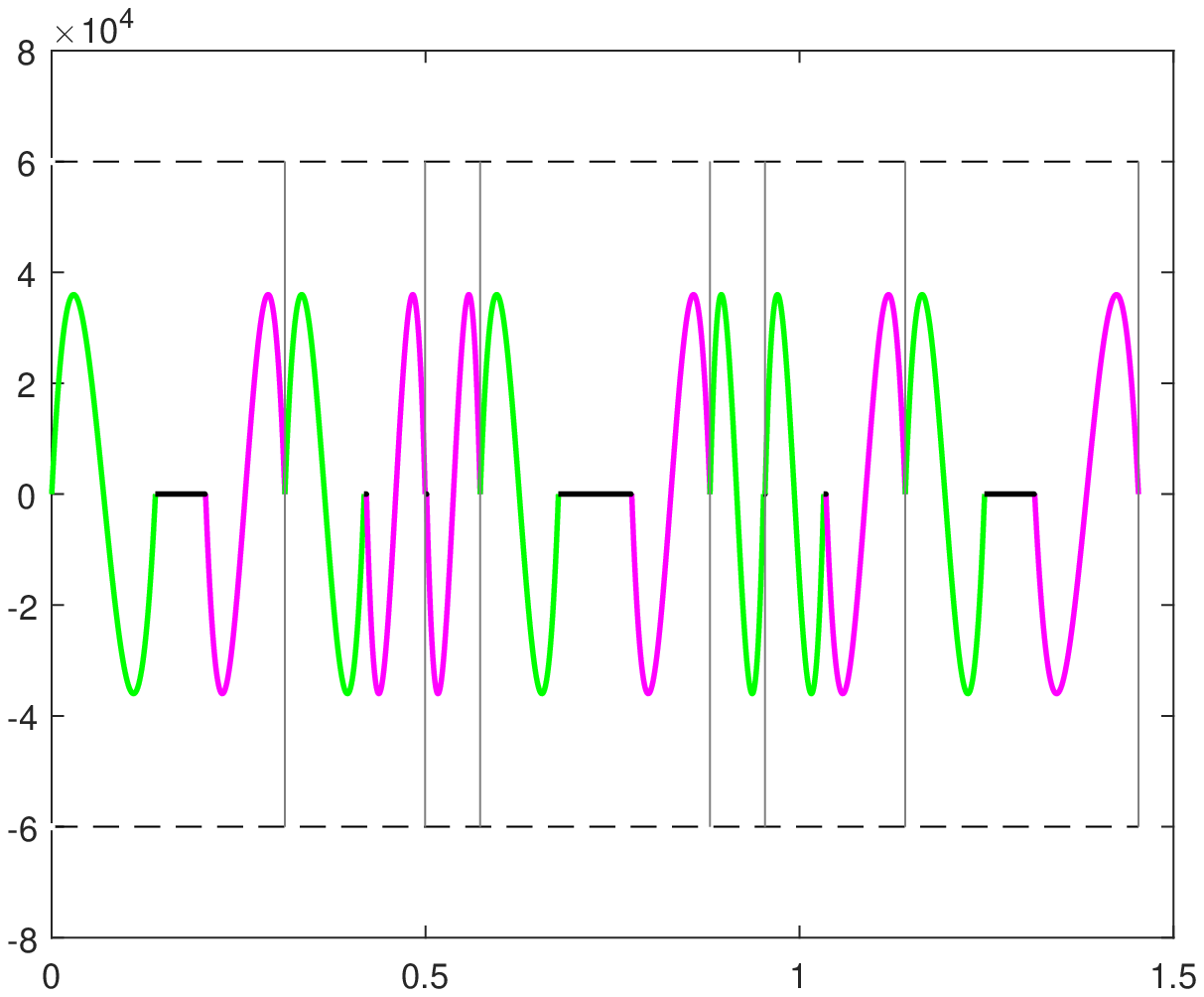}}
\caption{Feedrate scheduling results for the starfish curve obtained using the $R_1$ configuration. The chord error (a) and the feedrate profile (b) are shown together with the first (c) and second (d) derivatives of the feedrate function.}
\label{fig:star_feedrate}
\end{center}
\end{figure}

\begin{figure}[h]
\begin{center}
\subfigure[]{\includegraphics[scale=0.40,trim=0cm 0cm -1cm 0cm,clip]{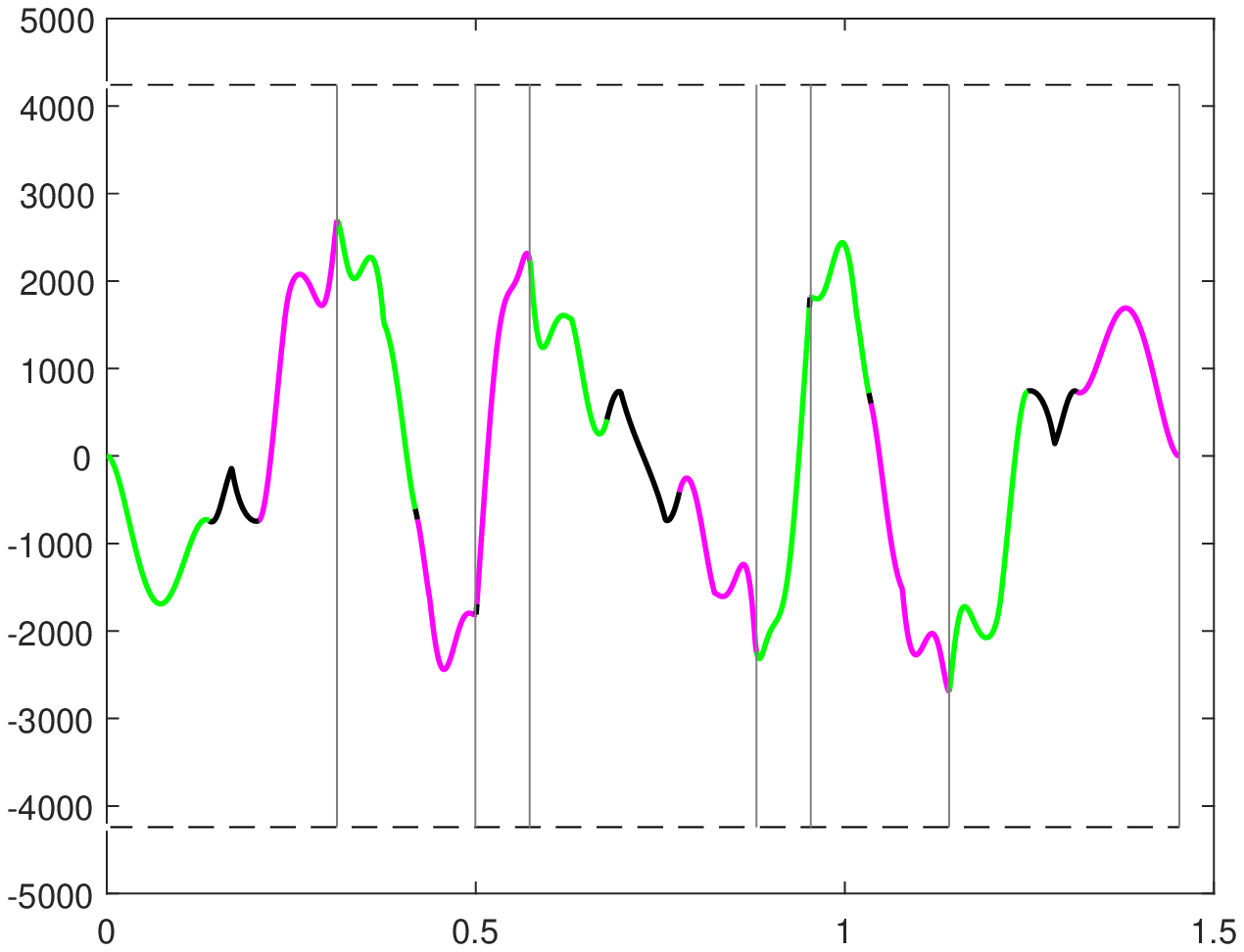}}\ \ \
\subfigure[]{\includegraphics[scale=0.40]{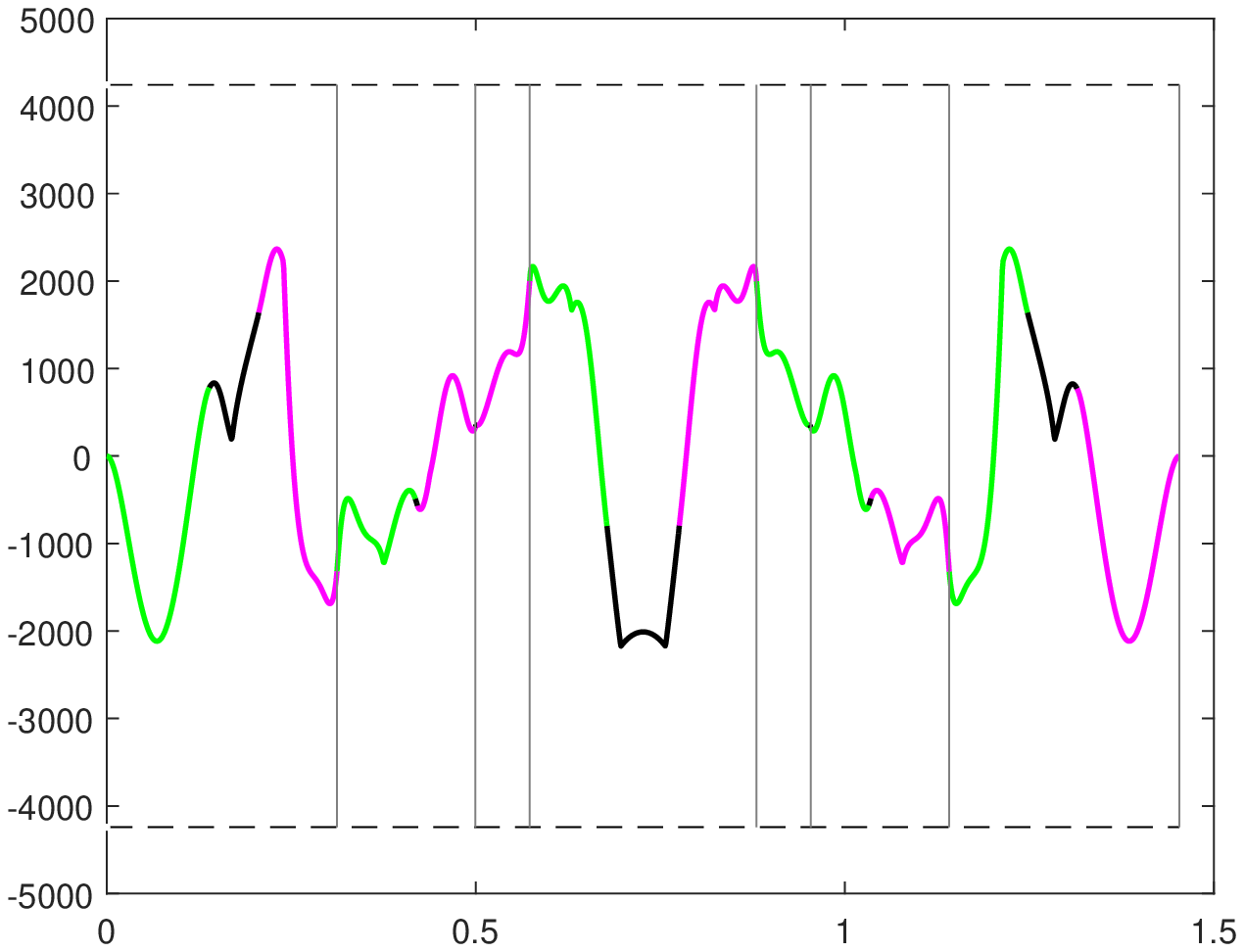}}\\
\subfigure[]{\includegraphics[scale=0.40,trim=0cm 0cm 1.3cm 0cm,clip]{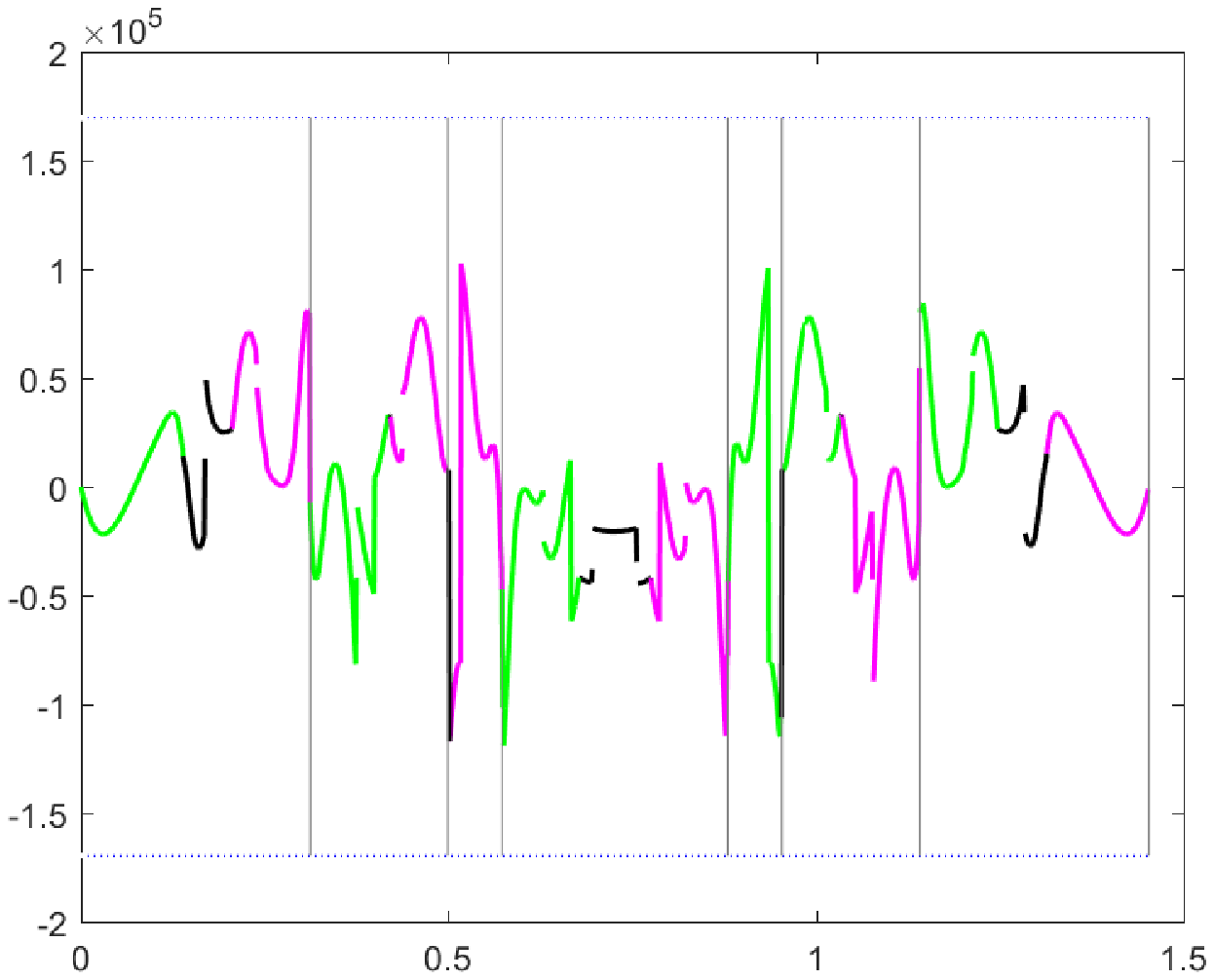}}\ \ \
\subfigure[]{\includegraphics[scale=0.40]{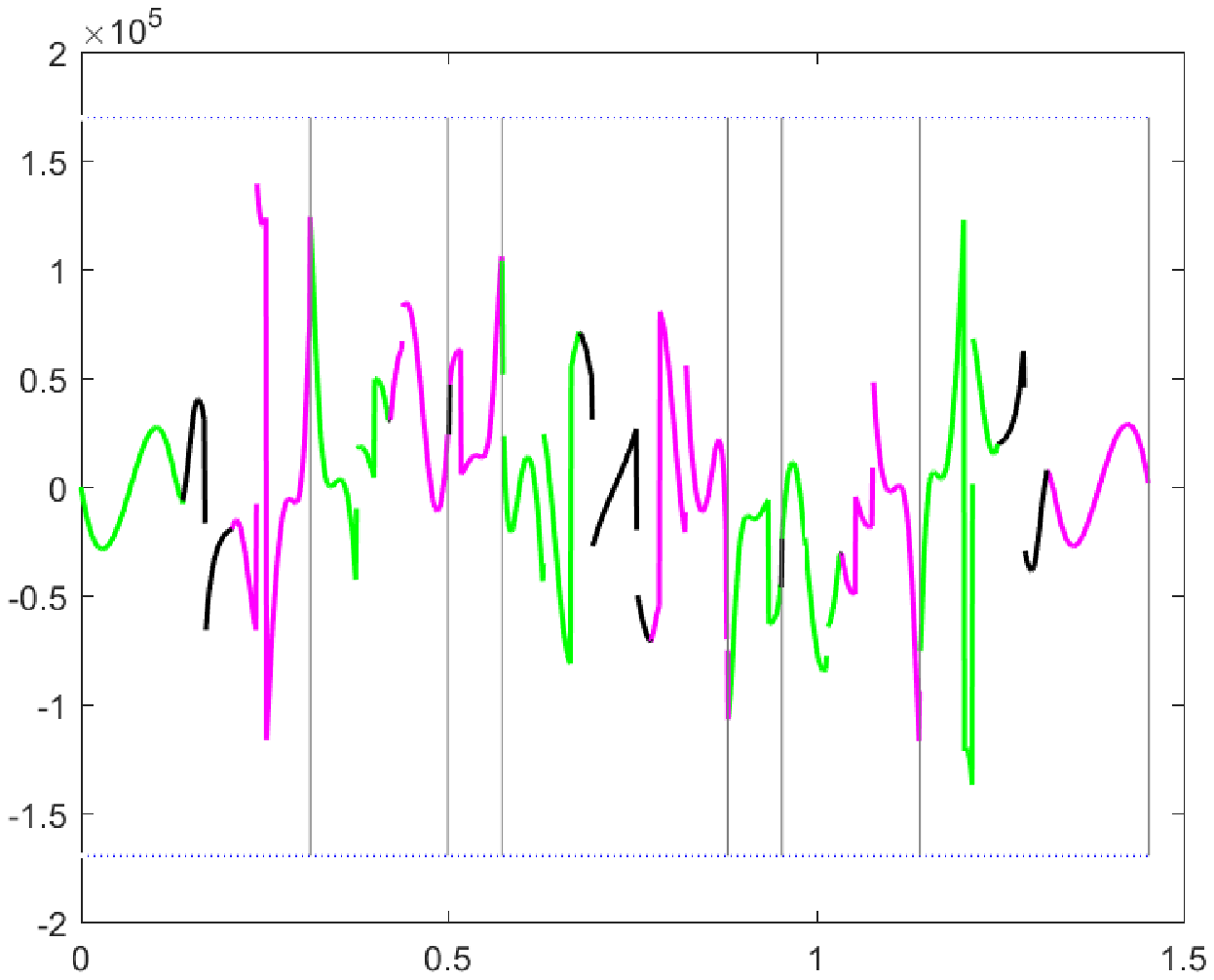}} 
\caption{Cartesian components of acceleration (a-b) and jerk (c-d)  obtained for the starfish curve with the $R_1$ configuration.}
\label{fig:star_cartesian}
\end{center}
\end{figure}

\begin{figure}[!h]
\begin{center}
\includegraphics[scale=0.33]{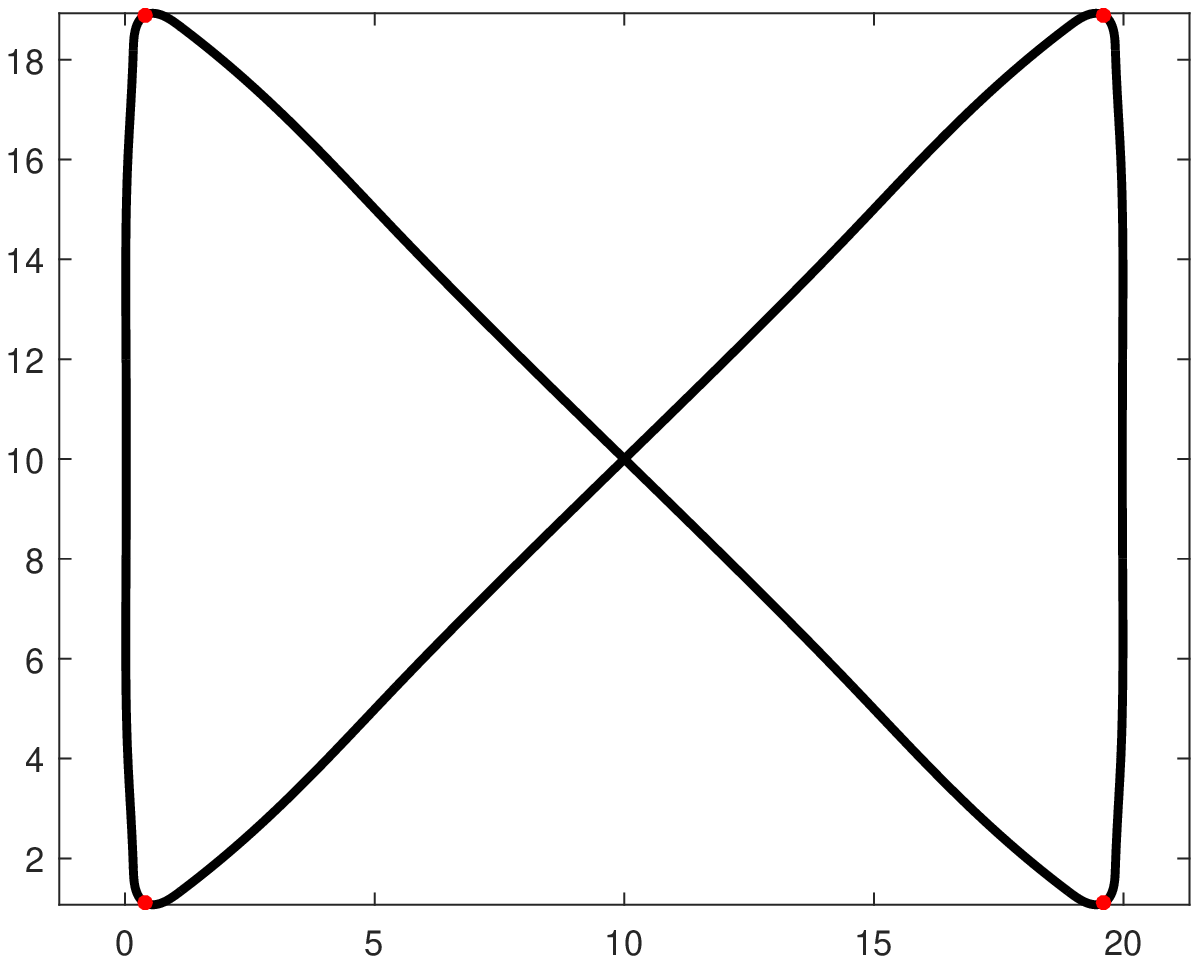}
\includegraphics[scale=0.33,trim=0cm 0.0cm 0cm 0cm,clip]{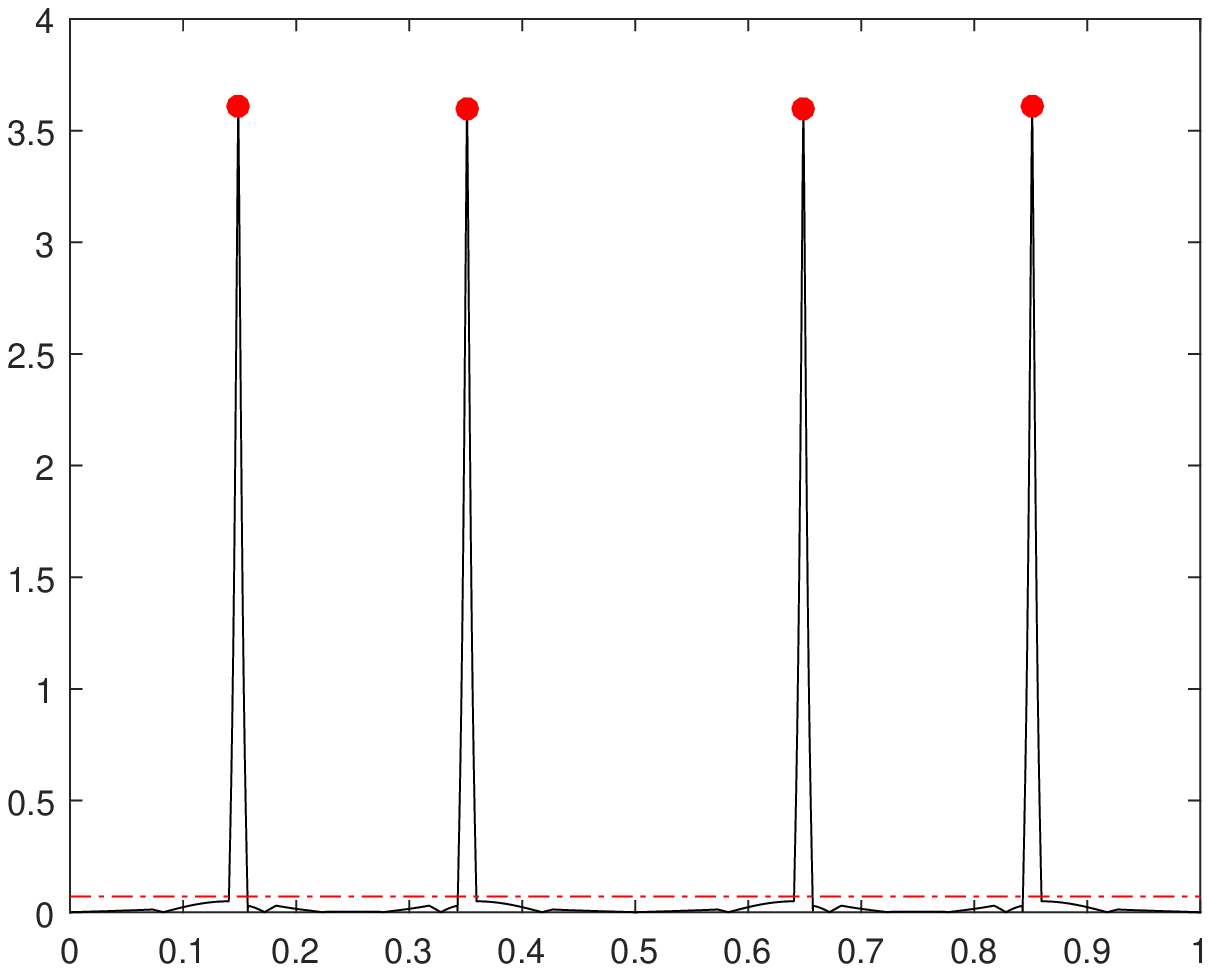}
\includegraphics[scale=0.33,trim=1.1cm 0.69cm -1cm 0cm,clip]{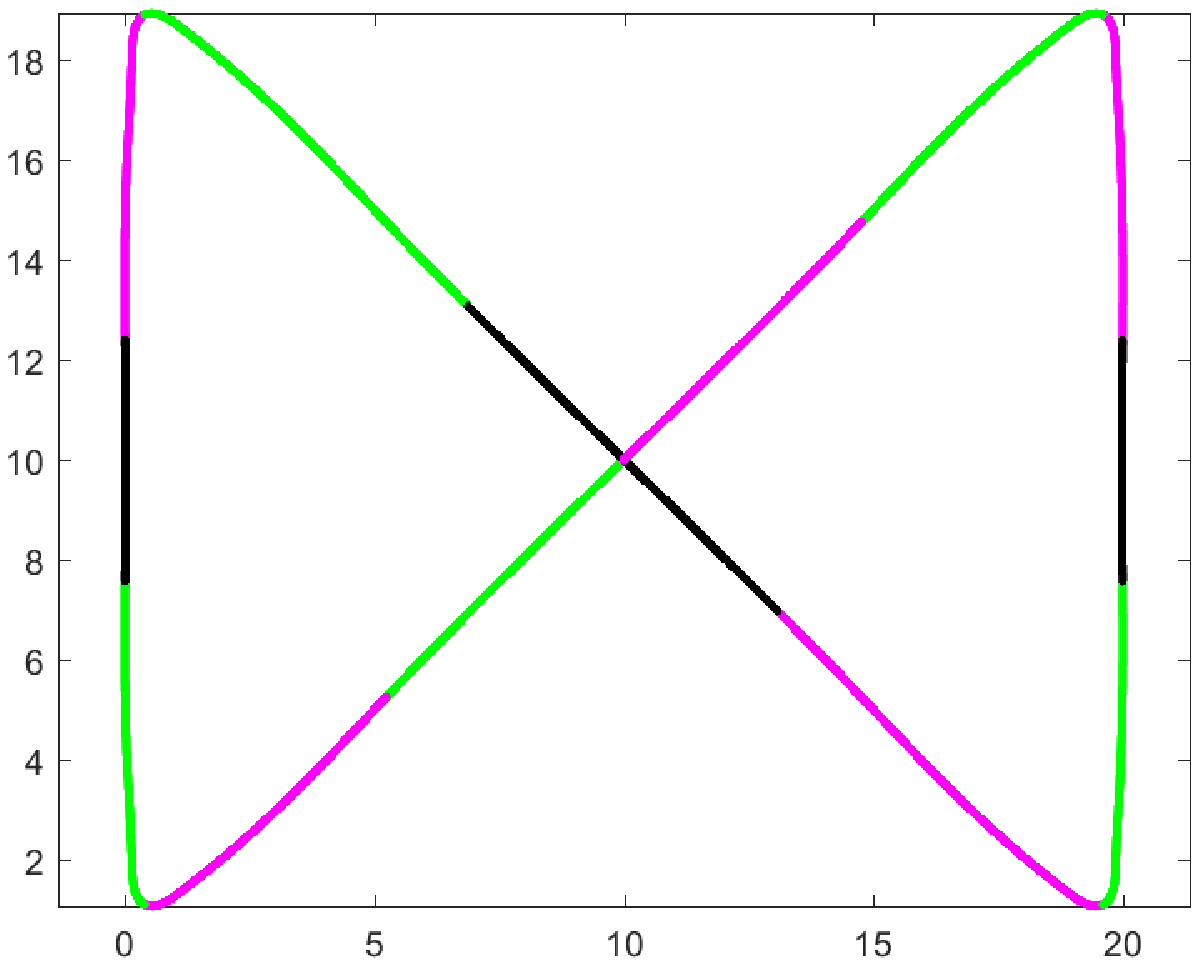}
\caption{Butterfly curve example: PH spline approximation (left), absolute curvature plot and critical curvature value (center, black and dashed red lines), and reference points (right). The critical points on the curve and the corresponding peaks on the
curvature plot are also shown (red dots).}
\label{fig:butterfly}
\end{center}
\end{figure}

\begin{figure}[h]
\begin{center}
\subfigure[]{\includegraphics[scale=0.40,trim=0cm 0cm -1cm 0cm,clip]{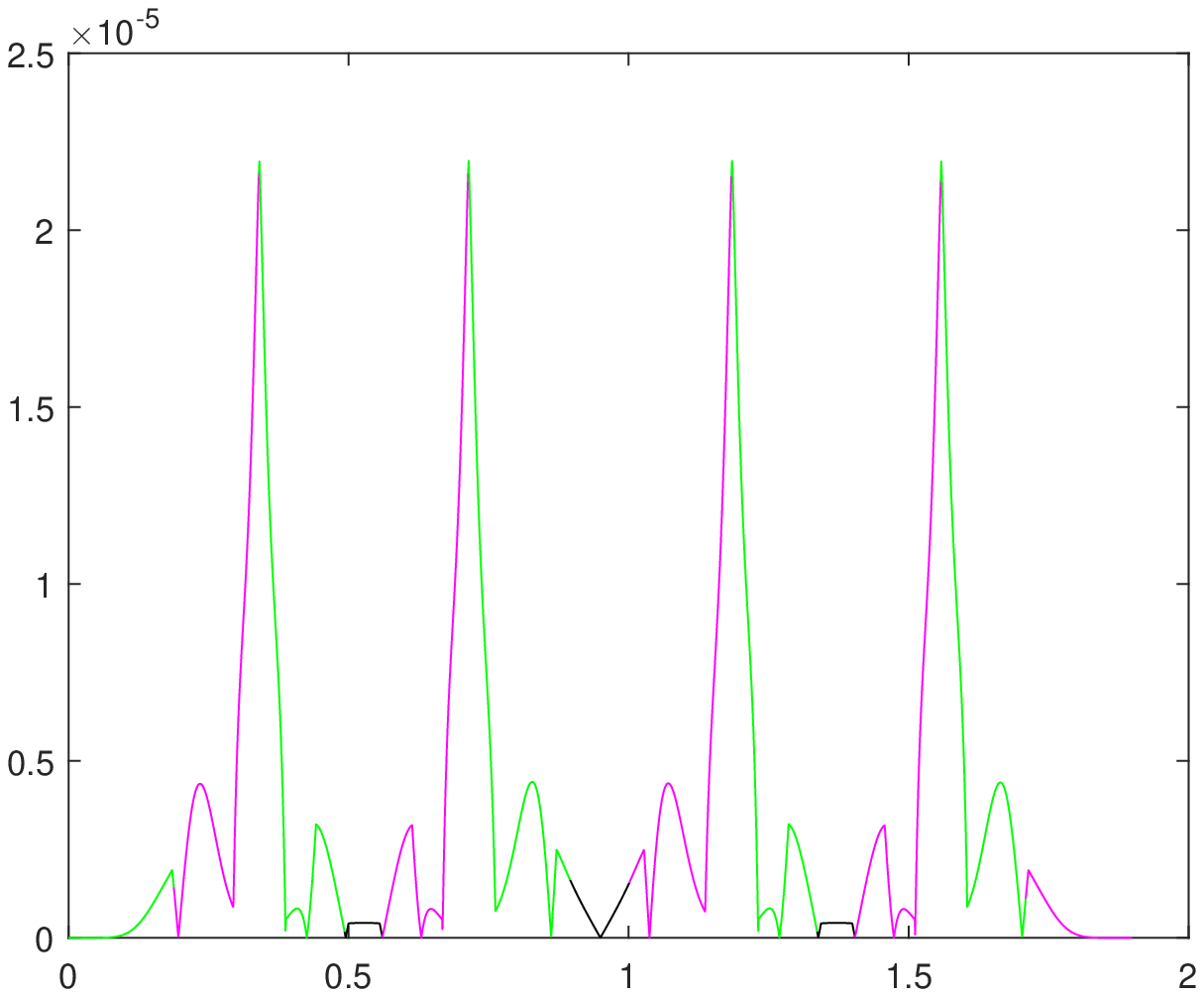}}\ \ \
\subfigure[]{\includegraphics[scale=0.40]{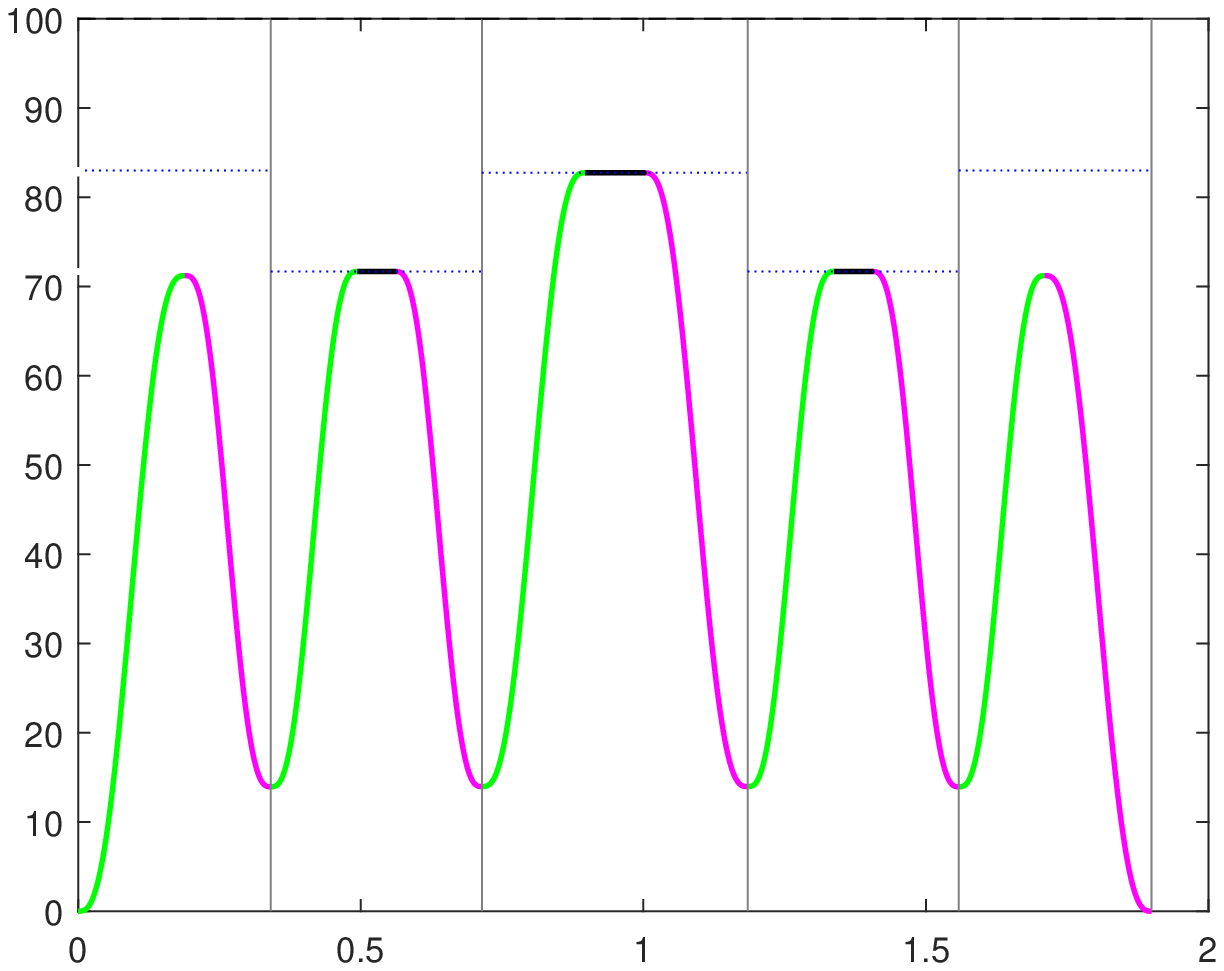}}\\
\subfigure[]{\includegraphics[scale=0.40,trim=0cm 0cm -1cm 0cm,clip]{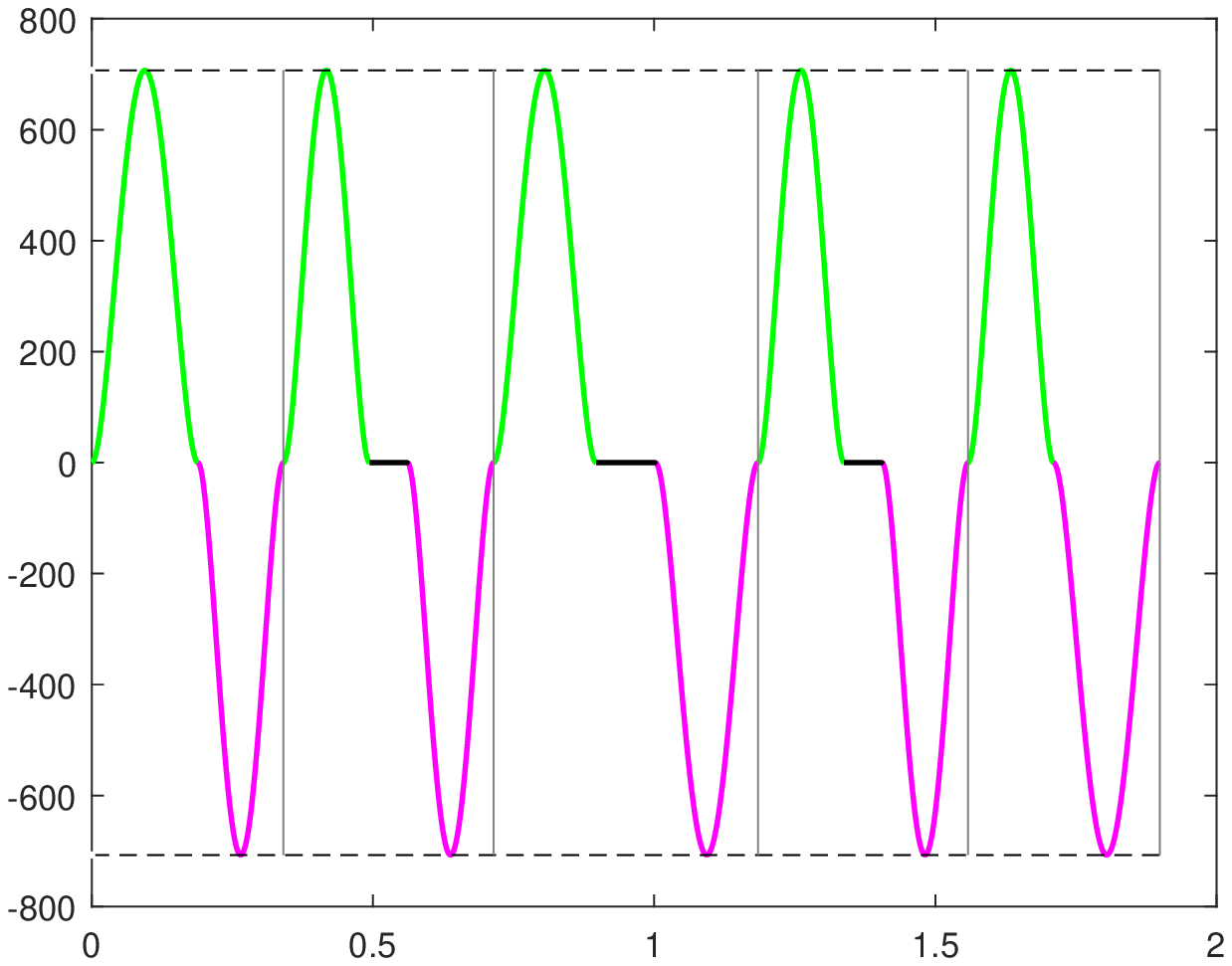}}\ \ \
\subfigure[]{\includegraphics[scale=0.40]{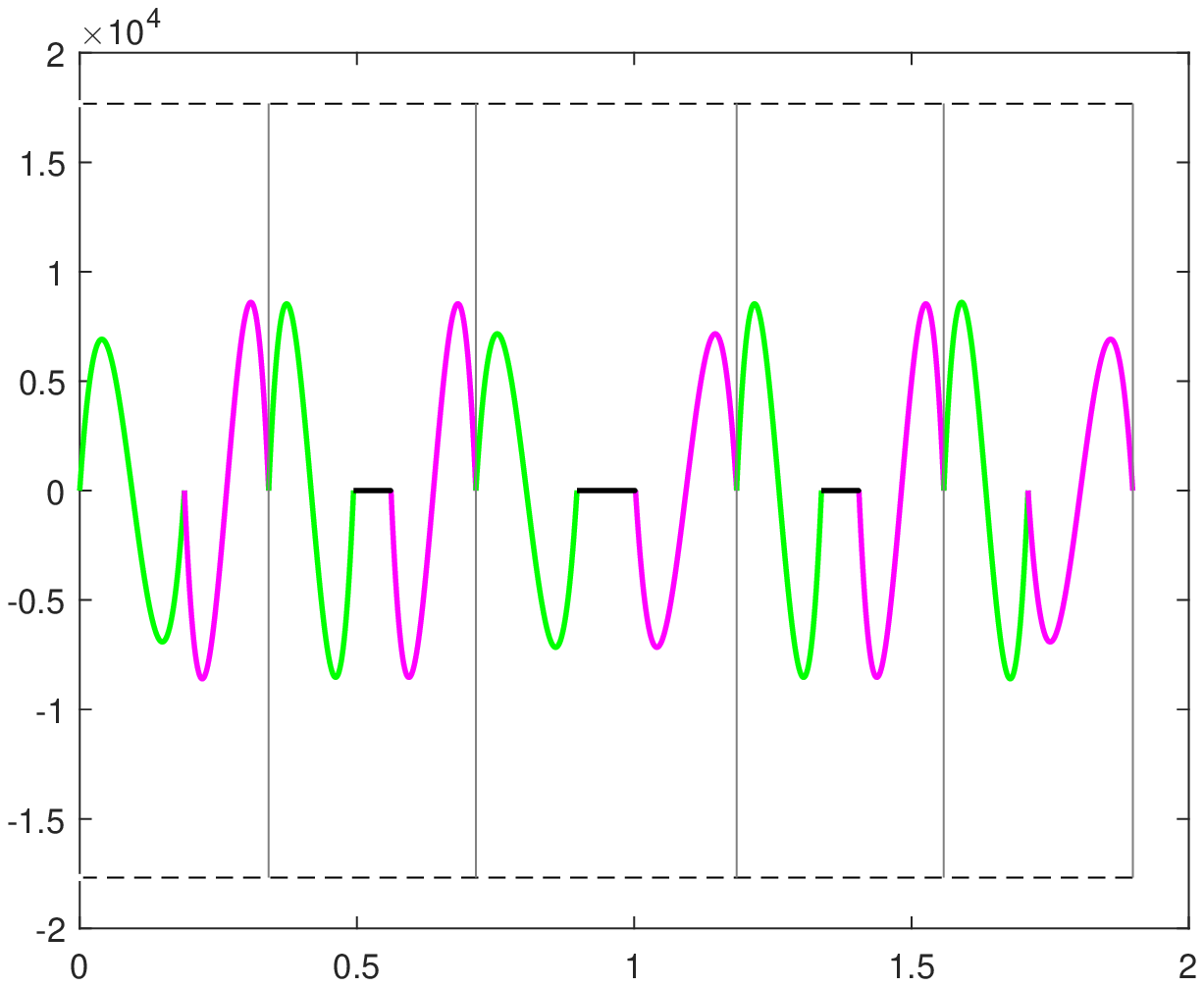}}
\caption{Feedrate scheduling results for the butterfly curve obtained using the $R_1$ configuration. The chord error (a) and the feedrate profile (b) are shown together with the first (c) and second (d) derivatives of the feedrate function.}
\label{fig:butterfly_feedrate}
\end{center}
\end{figure}

\begin{figure}[h]
\begin{center}
\subfigure[]{\includegraphics[scale=0.40,trim=0cm 0cm -1cm 0cm,clip]{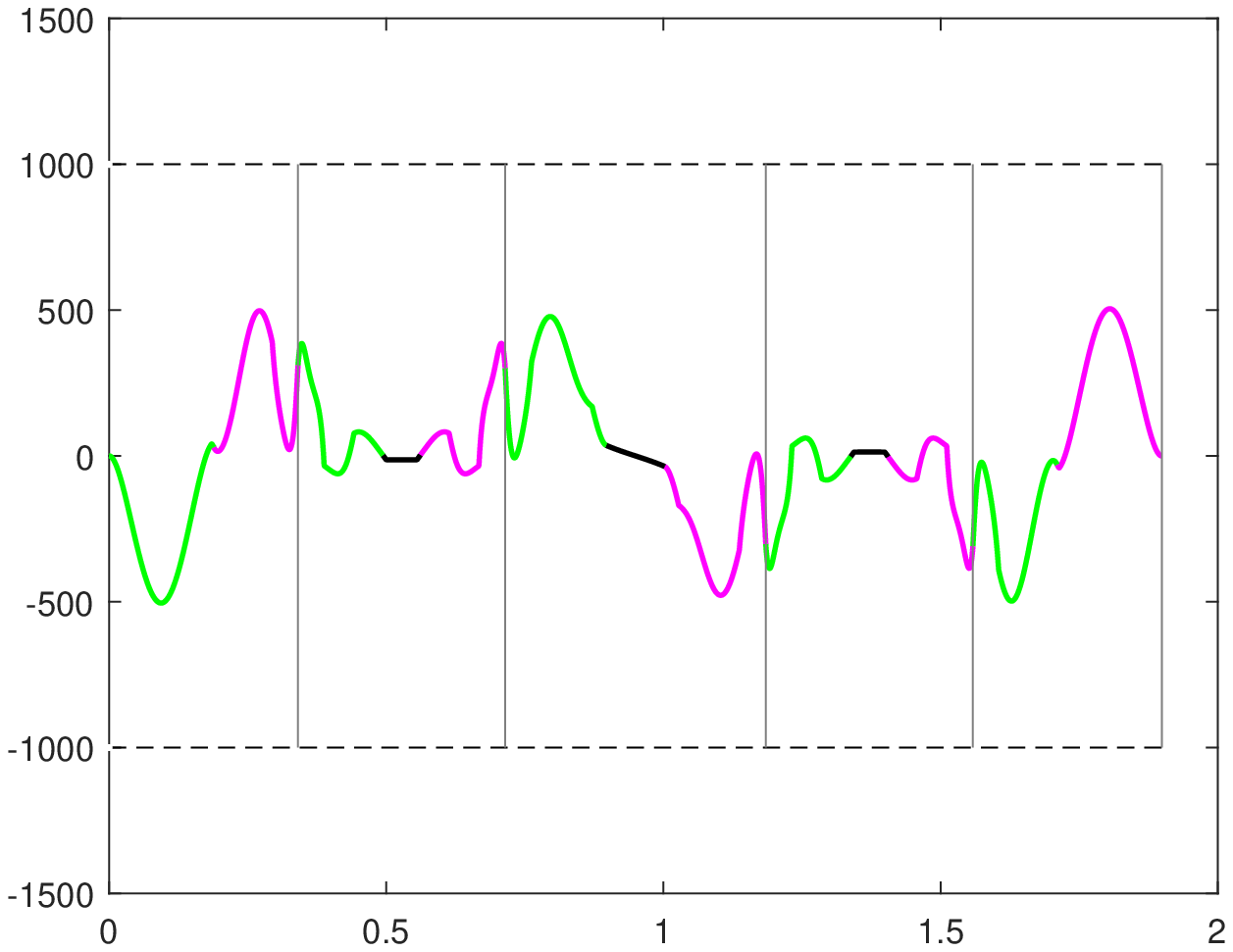}}\ \ \
\subfigure[]{\includegraphics[scale=0.40]{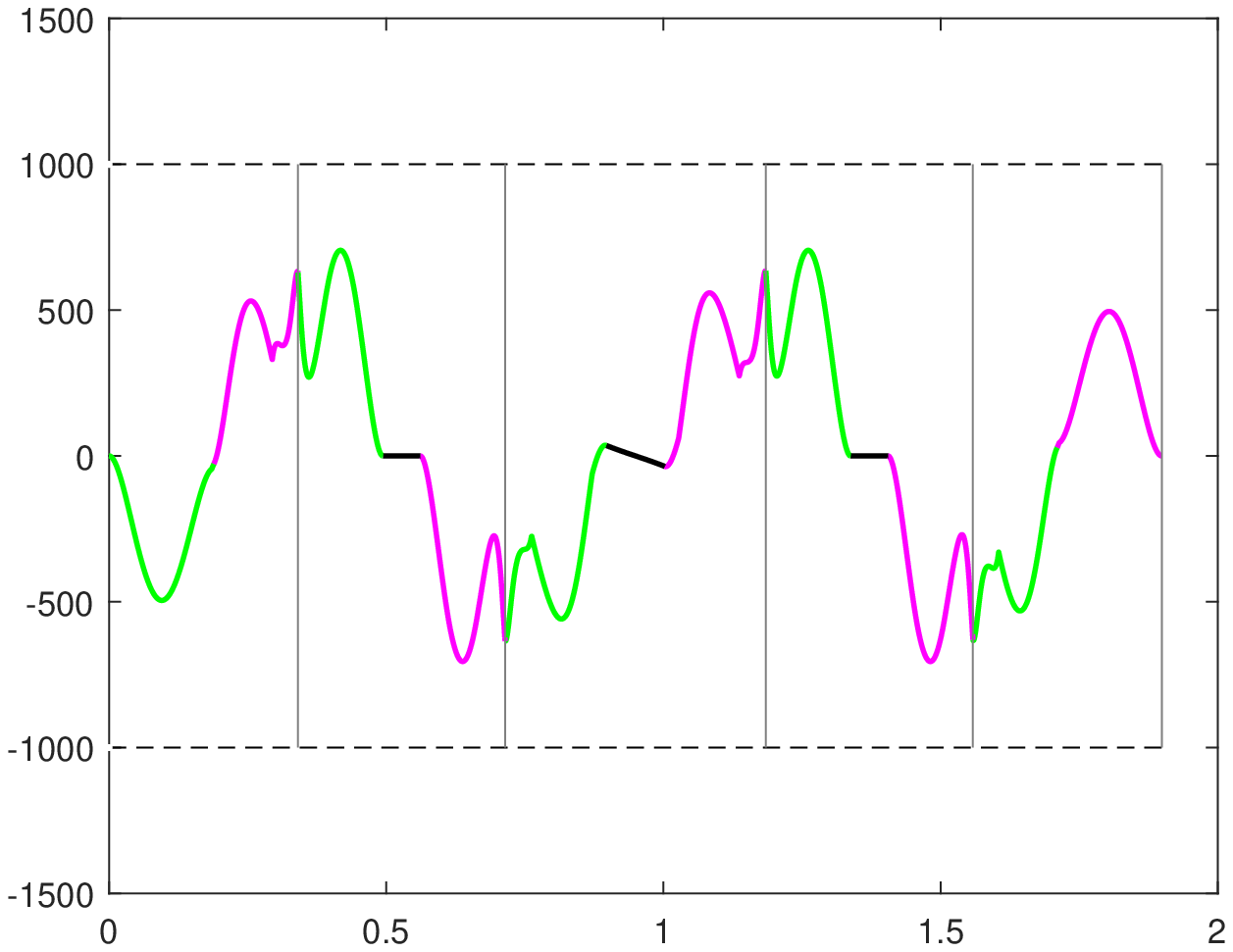}}\\
\subfigure[]{\includegraphics[scale=0.40,trim=0cm 0cm 1.3cm 0cm,clip]{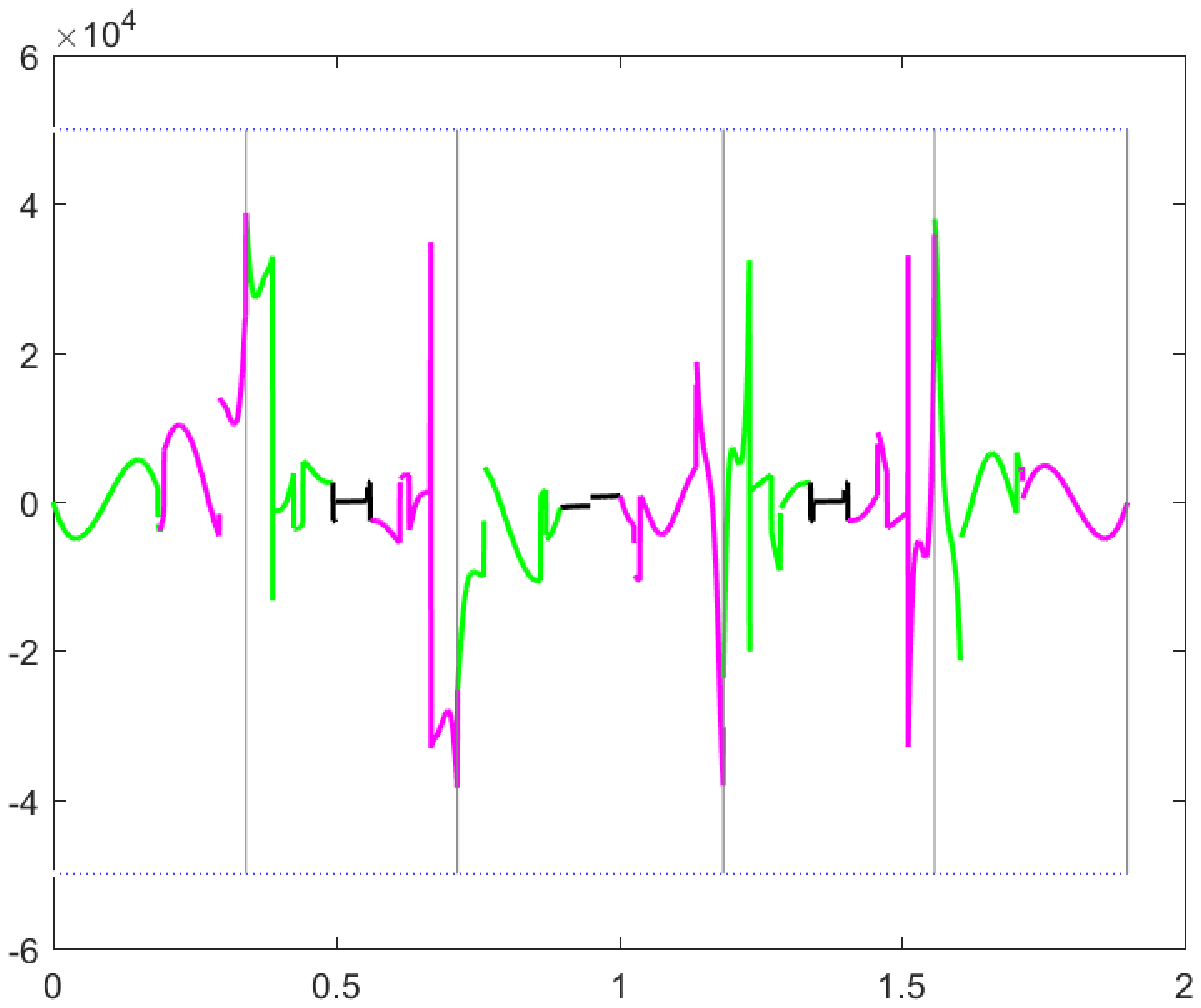}}\ \ \
\subfigure[]{\includegraphics[scale=0.40]{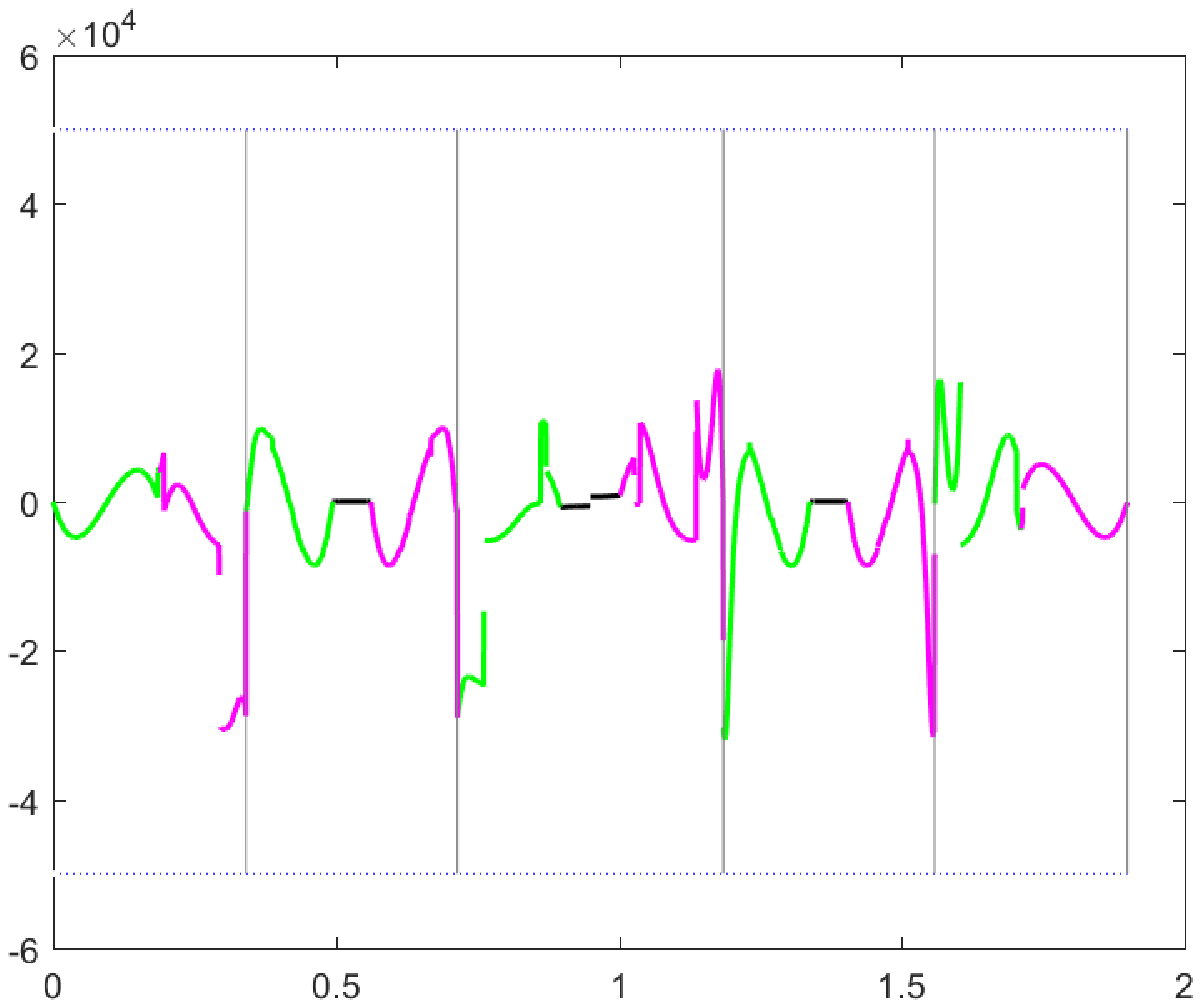}} 
\caption{Cartesian components of acceleration (a-b) and jerk (c-d)  obtained for the butterfly curve with the $R_1$ configuration.}
\label{fig:butterfly_cartesian}
\end{center}
\end{figure}

Figures \ref{fig:hat}, \ref{fig:star}, and \ref{fig:butterfly} show the PH spline approximation, the plot of the absolute value of the curvature, and the reference points obtained for the three test cases. The critical points and the corresponding peaks on the curvature plots are marked in red. The reference points associated to the acceleration, constant, and the deceleration phases are depicted in green, black, and magenta, respectively.
 Figures~\ref{fig:hat_feedrate}, \ref{fig:star_feedrate},
and \ref{fig:butterfly_feedrate} show the chord error and the resulting profiles of $v(t)$, $\dot{v}(t)$
and $\ddot{v}(t)$. The 
Cartesian components of acceleration and jerk are illustrated in Figures~\ref{fig:hat_cartesian},
\ref{fig:star_cartesian} and \ref{fig:butterfly_cartesian}. Again, acceleration, constant, and deceleration phases are indicated in green, black, and magenta, respectively.

For each example, the maximal feedrate $V_c^{(i)}$ is less
than $V_m$ in several curve blocks, due to the adaptive behaviour of our strategy described in Section
\ref{subsec:init_vm}, see Figures \ref{fig:hat_feedrate},
\ref{fig:star_feedrate} and \ref{fig:butterfly_feedrate}.
We can appreciate that the greater is the curvature on the $i$-th block, more
restrictive is the adjusted commanded feedrate. This behaviour is highlighted for
example on the second block of  the hat curve, where the curve presents an high
curvature peak and, consequently, the commanded feedrate drops to $115$.


Concerning the acceleration control, for all considered
examples the absolute value of the x-y acceleration components are always below
$A_m$, as shown in the first row of Figures \ref{fig:hat_cartesian},
\ref{fig:star_cartesian}, and \ref{fig:butterfly_cartesian}. In
particular, we can observe that on the hat curve the acceleration components
values  are rather below the prescribed threshold. In \cite{lee11} instead there exists one point where the threshold is not respected.
Furthermore the maximum chord error
$ D = 10^{-3}$mm is always fulfilled in our examples since the maximal chordal error is
$O(10^{-4})$, see Figures
\ref{fig:hat_feedrate}(a), \ref{fig:star_feedrate}(a), and
\ref{fig:butterfly_feedrate}(a).
The reference points produced by the PH interpolator with the sampling time
$\Delta t$ reported in Table \ref{tab:interp_param} are show on the right of Figures
\ref{fig:hat}, \ref{fig:star} and \ref{fig:butterfly}.

Finally, for all the three considered tests, Figures~\ref{fig:hat_feedrate},
\ref{fig:star_feedrate} and \ref{fig:butterfly_feedrate}(d) and Figures
\ref{fig:hat_cartesian}, \ref{fig:star_cartesian} and
\ref{fig:butterfly_cartesian} show the obtained profiles of $\ddot
v$ and the x-y jerk components, respectively. The
desired control on $\ddot v$ is ensured by our motion in all the cases. The x-y jerk components are also everywhere suitably controlled for each test.
 

\section{Closure}\label{sec:clo}

We introduced a configurable $C^2$ feedrate scheduler whose more restrictive formulation ensures full kinematic control.  Faster motions are obtained by using relaxed configurations, which often still preserve the desired motion control. The  developed scheduler has been applied to $C^2$ PH spline curves to produce acceleration continuous motions.  The selection of  paths described in terms of PH spline curves is of great interest since it  avoids numerical approximations of fundamental geometrical quantities and fosters the accuracy and robustness of the overall scheme.

The focus of this paper has been on the feedrate profile for planar (PH) paths. Although the spatial case requires a separate study, the procedure should readily extends to the three-dimensional setting. In order to obtain a jerk continuous motion  on a PH spline path, a $C^3$ Hermite interpolation scheme with PH curves (of higher degree) should be investigated.

\section*{Acknowledgments}

This work has been partially supported by the MIUR ``Futuro in Ricerca'' project DREAMS (RBFR13FBI3), and by Istituto Nazionale di Alta Matematica (INdAM) through  Finanziamenti Premiali SUNRISE.

\bibliographystyle{abbrv}
\bibliography{biblio}

\end{document}